\documentclass[righttag,12pt]{amsart}
\usepackage{amssymb}
\usepackage{euscript}

\usepackage[all]{xy}
\xyoption{arc}

\newtheorem{prop}{Proposition}[section]
\newtheorem{theorem}[prop]{Theorem}
\newtheorem{cor}[prop]{Corollary}
\newtheorem{lemma}[prop]{Lemma}
\newtheorem{example}[prop]{Example}
\newtheorem{step}{Step}
 
\numberwithin{equation}{section}

\theoremstyle{definition}
\newtheorem{fig}{\hspace{65mm}{\sc Figure}}[section]
\newtheorem{defn}[prop]{Definition}
\newtheorem{ack}{Acknowledgments}

\theoremstyle{remark}
\newtheorem{rem}[prop]{Remark}
\newcommand{\bbN}{{\mathbb{N}}}
 \newcommand{\bbR}{{\mathbb{R}}}
\newcommand{\bbZ}{{\mathbb{Z}}}

\newcommand{\al}{\alpha}
\newcommand{\be}{\beta}

\newcommand{\de}{\delta}
\newcommand{\e}{\varepsilon}

\newcommand{\ro}{\varrho}

 \newcommand{\Om}{\Omega}

\newcommand{\card}{\operatorname{card}}
\newcommand{\Bd}{\operatorname{Bd}}

\newcommand{\lb}{\label}

\newcommand{\lra}{\longrightarrow}

\newcommand{\DEF}{\buildrel {\mbox{\small def}}\over =}

\begin{document}

\title{On the structure of level sets of uniform and Lipschitz
quotient mappings from ${\mathbb{R}}^n$ to ${\mathbb{R}}$}

\author{Beata Randrianantoanina$^*$}\thanks{$^*$Participant, NSF Workshop
in Linear Analysis and Probability, Texas A\&M University}

\address{Department of Mathematics and Statistics \\ Miami University
\\Oxford, OH 45056}

\email{randrib@muohio.edu}


\begin{abstract}
We study two questions posed by Johnson, Lindenstrauss, Preiss, and Schechtman,
concerning the structure of level sets of uniform and Lipschitz quotient mappings from
$\bbR^n \to \bbR$.
We show that if $f: {\mathbb{R}}^n \longrightarrow {\mathbb{R}}$, $n \geq 2$,
is a uniform quotient
mapping then for
every $t \in {\mathbb{R}}$, $f^{-1}(t)$ has a bounded number of components, each component of
$f^{-1}(t)$ separates ${\mathbb{R}}^n$ and the upper bound of the number of components depends
only on $n$ and the moduli of co-uniform and uniform continuity of $f$.

Next we obtain  a characterization
of the form of any
  closed, hereditarily locally connected, locally compact, connected set with no end
 points and containing no simple closed curve, and we apply it to
 describe the structure of level sets of co-Lipschitz uniformly continuous mappings $f: {\mathbb{R}}^2
\longrightarrow  {\mathbb{R}}$. We prove that  all level sets of any co-Lipschitz uniformly continuous
mapping
from $ {\mathbb{R}}^2$ to
${\mathbb{R}}$ are locally connected, and we show that
for every pair of a constant $c>0$ and a function $\Om(\cdot)$ with $\lim_{r\to 0}\Om(r)=0$,
there exists a natural number $M=M(c,\Om)$, so
that for every co-Lipschitz uniformly continuous map $f:\bbR^{2} \lra \bbR$ with a co-Lipschitz
constant $c$
and a modulus of uniform continuity $\Om$,
 there exists a natural number $n(f)\leq M$ and a finite set $T_f \subset {\mathbb{R}}$
with $\card(T_f) \leq  n (f)-1$ so that for all $t \in {\mathbb{R}} \setminus T_f$,
$f^{-1}(t)$ has exactly $n(f)$ components, ${\mathbb{R}}^2 \setminus f^{-1}(t)$ has exactly
$n(f)+ 1$ components and each component of $f^{-1}(t)$ is homeomorphic with the real
line and separates the plane into exactly $2$ components.  The number and form of
components of $f^{-1}(s)$ for $s \in T_f$ are also described -- they have a
finite graph structure.

We give an example of a uniform quotient map from ${\mathbb{R}}^2$ to ${\mathbb{R}}$
 which has
non-locally connected level sets.
\end{abstract}

 \maketitle

\section{Introduction}
Let $X,Y$ be metric spaces.  A mapping $f: X \longrightarrow Y$ is said to be a {\it
uniform quotient mapping} if there exist functions $\omega, \Omega: \bbR_+ \to \bbR_+$
with $\omega(r)>0$ for all $r>0$ and $\lim_{r \to 0} \Omega(r)=0$ so that for
all $x \in X$ and all $r>0$:
\begin{equation}\lb{uqm}
B(f(x),\omega(r))\subset f (B(x,r))\subset B(f(x), \Omega(r)),
\end{equation}
where $B(x,r)$ denotes the open ball with center $x$ and radius $r$.

Notice that the right hand inclusion means that $f$ is uniformly continuous.  The
mapping $f$ is called {\it co-uniformly continuous} if the left hand inclusion in
\eqref{uqm} is satisfied.  There is no restriction in assuming that the functions $\omega$
and $\Omega$ are continuous and increasing.  They are called {\it moduli of co-uniform
and uniform continuity of $f$}, respectively.  If the functions $\omega$ and $\Omega$ are
linear, i.e. if there exist constants $c, L>0$ so that for all $x \in X$ and all
$r>0$:
\begin{equation}\lb{Lqm}
B(f(x),c r) \subset f(B(x,r))\subset B(f(x),L r),
\end{equation}
then $f$ is called a {\it Lipschitz quotient mapping}.  Clearly the right hand
inclusion in \eqref{Lqm} means that $f$ is a Lipschitz mapping.  If $f$ satisfies the
left hand inclusion of \eqref{Lqm}, $f$ is called a {\it co-Lipschitz mapping}.
Constants $c$ and $L$ are called {\it co-Lipschitz and Lipschitz constants of $f$},
respectively.  The study of uniform and Lipschitz quotient mappings was initiated in
\cite{BJLPS}, see also \cite{BL} for the comprehensive introduction of the subject.
The structure of Lipschitz and uniform quotient mappings $f: X \longrightarrow Y$, when $X$
and $Y$ are finite dimensional was studied by Johnson, Lindenstrauss, Preiss and
Schechtman in \cite {JLPS00}.  They obtained   most complete results for the case of
$X = Y = \bbR^2$. For $f:\bbR^2 \lra \bbR^2$ they proved, in particular, that if $f$ is
a uniformly continuous and co-Lipschitz, e.g. if
$f$ is a Lipschitz quotient mapping, then for every $t \in \bbR^2$, $f^{-1}(t)$ is a
finite set of points in $\bbR^2$ and $f=P \circ h$ where $h$ is a homeomorphism of the
plane and $P$ is a complex polynomial (see also Remark~\ref{rem-a} below).
The question whether level sets of $f^{-1}(t)$
of a Lipschitz quotient map $f: \bbR^n \lra \bbR^n$ are discrete, is open for all $n
>2$.

In \cite{JLPS00}, the authors also study the structure of level sets $f^{-1}(t)$ of
uniform and Lipschitz maps $f: \bbR^n \lra \bbR$.  They showed, among others,
 the following results:

\begin{theorem} \lb{finitecompl} \cite[Proposition~5.1]{JLPS00}
Let $f:\bbR^n \to \bbR$ be a uniform quotient mapping satisfying \eqref{uqm}.  Then
for each $t \in \bbR$ the number of components of $\bbR^n\setminus f^{-1}(t)$ is
finite and bounded by a function of $n$, $ \omega (\cdot)$ and $\Omega (\cdot)$ only.
\end{theorem}

\begin{theorem}\cite[Proposition~5.4]{JLPS00}\lb{sep}
Let $f: \bbR^2 \lra \bbR$ be a Lipschitz quotient mapping.  Then, for each $t \in
\bbR$, each component of $f^{-1}(t)$ is unbounded and  separates the plane.
\end{theorem}

\begin{theorem}\cite[Corollary~5.5]{JLPS00}\lb{fin}
Let $f: \bbR^2 \lra \bbR$ be a Lipschitz quotient mapping.  Then, for each $t \in \bbR,
f^{-1}(t)$ has a bounded number of components.  The upper bound of the number of
components depends only on the Lipschitz and co-Lipschitz constants of $f$.
\end{theorem}

They also asked the following two questions:
\begin{itemize}
\item[(Q1)] Can one weaken the assumption of Lipschitz quotient to uniform quotient mappings
in Theorems~\ref{sep} and \ref{fin}?
\item[(Q2)] To what extent is the number of components of $f^{-1}(t)$ or of
 $\bbR^2 \setminus f^{-1}(t)$ independent of $t$?  Are these numbers constant
 after excluding finitely many values of $t$?
\end{itemize}

Question (Q2) is motivated by the following two examples of Lipschitz quotient
mappings from $\bbR^2$ to $\bbR$.  In both cases the mapping $f$ is the $\ell_1$ distance from
the solid lines multiplied, in each component of the complement of the solid lines, by
the sign indicated.

\vspace{2mm}

 \hspace{2mm} \xy <1.2cm,0cm>:
 (-1.5,0) ;
 (1.5,0)**\dir{-},
(-1.5,0) ;
 (-3,0)**\dir{-},
 (1.5,0) ;
 (3,0)**\dir{-},
 (2,1.5) ;
 (1.3,.8)**@{-},
 (2,-1.5) ;
 (1.3,-.8)**@{-},
 (0.5,0) ;
(1.3,-.8)**@{-},
(0.5,0) ;
(1.3,.8)**@{-},
(2,1.5) ;
 (1.3,.8)**@{-},
 (-2,-1.5) ;
 (-1.3,-.8)**@{-},
 (-0.5,0) ;
(-1.3,-.8)**@{-},
(-0.5,0) ;
(-1.3,.8)**@{-},
(-2,1.5) ;
 (-1.3,.8)**@{-},
 (9,-1.5) ;
 (7.5,0)**@{-},
(7.5,0) ;
(9,1.5)**@{-},
 (5,-1.5) ;
 (6.5,0)**@{-},
(6.5,0) ;
(5,1.5) **@{-},
(-2,.6)*+{-},
(0,.6)*+{+},
(2,.6)*+{-},
(-2,-.6)*+{+},
(0,-.6)*+{-},
(2,-.6)*+{+},
(5.5,0)*+{-},
(7,0)*+{+},
(8.5,0)*+{-},
\endxy

\begin{fig}  \label{picture}
\end{fig}


Here $f^{-1}(0)$ has one component in the first example and two in the second, and
$\bbR^2 \setminus f^{-1}(0)$ has six and three components, respectively.  The authors
of \cite{JLPS00} note that it is easy to draw examples with an arbitrary finite number
of components of  $f^{-1}(0)$. Thus question (Q2) is essentially asking whether all
Lipschitz quotient maps $\bbR^2$ to $\bbR$ have the form similar to the  examples illustrated
in Figure~\ref{picture}.

This paper is devoted to the study of questions (Q1) and (Q2).  We answer both of them
affirmatively.
First, in Section~\ref{sec-sep}, we obtain generalizations of Theorems~\ref{fin} and
\ref{sep} for uniform quotient mappings from $\bbR^n$ to $\bbR$ for any $n \geq 2$
(Theorem~\ref{finite}) and Corollary~\ref{sepplane}, respectively). Our results follow
from Theorem~\ref{finitecompl} through general topological arguments based on the
Phragmen-Brower theorem and the theory of separation in $\bbR^n$.

Next we study question (Q2).  We obtain not only information about the number of
components of $f^{-1}(t)$ and of $\bbR^2 \setminus f^{-1}(t)$ for Lipschitz quotient maps
$f: \bbR^2 \lra \bbR$, but we give the full characterization of both the number and the
form of each component of any level set $f^{-1}(t)$ of co-Lipschitz uniformly continuous
mappings $f:\bbR^{2} \lra \bbR$ (Theorem~\ref{main}).  We show that
for every pair of a constant $c>0$ and a function $\Om(\cdot)$ with $\lim_{r\to 0}\Om(r)=0$,
there exists a natural number $M=M(c,\Om)$, so
that for every co-Lipschitz uniformly continuous map $f:\bbR^{2} \lra \bbR$ with a co-Lipschitz
constant $c$
and a modulus of uniform continuity $\Om$,
 there exists a natural number $n(f)\leq M$ and a finite set $T_f \subset \bbR$
with $\card(T_f) \leq   n (f)-1$ so that for all $t \in \bbR \setminus T_f$,
$f^{-1}(t)$ has exactly $n(f)$ components, $\bbR^2 \setminus f^{-1}(t)$ has exactly
$n(f)+ 1$ components and each component of $f^{-1}(t)$ is homeomorphic with the real
line and separates the plane into exactly $2$ components.  The number and form of
components of $f^{-1}(s)$ for $s \in T_f$ is also described -- these components have a
finite graph structure (for precise formulation see Theorems~\ref{main},
\ref{graph} and Remark~\ref{rem-n}).

Thus we do confirm that co-Lipschitz uniformly continuous mappings from $\bbR^2$ to $\bbR$ have a
form analogous to the examples presented on Figure~\ref{picture}.  Moreover, we prove
that, as on Figure~\ref{picture}, no level set $f^{-1}(t)$ can contain parallel lines,
but the distance between unbounded components of $f^{-1}(t)\setminus B(0,R)$ has to increase to
infinity as $R$ increases to infinity, cf. Figure~\ref{pic2}
(Proposition~\ref{faraway}).

\vspace{3mm}

  \hspace{2mm} \xy <1.2cm,0cm>:
 (-1.5,0) ;
 (1.5,0)**\dir{.},
(-1.5,0) ;
 (-3,0)**\dir{-},
 (1.5,0) ;
 (3,0)**\dir{-},
 (2,1.5) ;
 (1.3,.8)**@{-},
 (2,-1.5) ;
 (1.3,-.8)**@{-},
 (0.5,0) ;
(1.3,-.8)**@{.},
(0.5,0) ;
(1.3,.8)**@{.},
(2,1.5) ;
 (1.3,.8)**@{-},
 (-2,-1.5) ;
 (-1.3,-.8)**@{-},
 (-0.5,0) ;
(-1.3,-.8)**@{.},
(-0.5,0) ;
(-1.3,.8)**@{.},
(-2,1.5) ;
 (-1.3,.8)**@{-},
(0,0)*\ellipse<18mm>{--},
 (9,1.5) ;
 (8.3,.8)**@{-},
 (9,-1.5) ;
 (8.3,-.8)**@{-},
 (7.5,0) ;
(8.3,-.8)**@{.},
(7.5,0) ;
(8.3,.8)**@{.},
(9,1.5) ;
 (8.3,.8)**@{-},
 (5,-1.5) ;
 (5.7,-.8)**@{-},
 (6.5,0) ;
(5.7,-.8)**@{.},
(6.5,0) ;
(5.7,.8)**@{.},
(5,1.5) ;
 (5.7,.8)**@{-},
(3.5,0)*\ellipse<18mm>{--},
\endxy
\begin{fig}  \label{pic2}
\end{fig}

Our method of proof of Theorem~\ref{main} depends on a careful analysis of topological
properties of level sets  $f^{-1}(t)$, their end points and their structure at
infinity.  The crucial property that we use in a very essential way is the fact that
level sets $f^{-1}(t)$ are locally connected when $f$ is a co-Lipschitz uniformly
continuous map from
$\bbR^2$ to $\bbR$ (Proposition~\ref{r2}).

  Our first  characterization of the structure of level sets is in fact a characterization
of the form of any
  closed, hereditarily locally connected, locally compact, connected set with no end
 points and containing no simple closed curve.
We present a new self-contained proof of this characterization. Similar characterizations
for dendrites and for sets whose
every point is a cut point have been previously obtained by Shimrat \cite{Sh59},
Stone
\cite{S68} and  Nadler \cite{N93} (see Theorem~\ref{graph}
and Remarks~\ref{rem-d}, \ref{rem-tree}).

We do not know whether level sets of
co-Lipschitz uniformly continuous maps or of Lipschitz quotient maps from
$\bbR^n$ to $\bbR$ are locally connected when $n > 2$.  If
one looks for a counter-example, the most natural map to check would be the Lipschitz
quotient map $f:\bbR^3 \lra \bbR^2$ constructed by Cs\"ornyei \cite{C01}, whose level
set $f^{-1}(0)$ is very large and complicated.  It turns out, however, that for this
map and also for its both coordinate maps, which go from $\bbR^3$ to $\bbR$, all level
sets are locally connected.

However we do know that there exist
 uniform quotient maps from $\bbR^2$ to $\bbR$ with non-locally
connected level sets (see Example~\ref{exnlc}).

The local connectedness of level sets $f^{-1}(t)$ of a co-Lipschitz uniformly continuous map
$f:\bbR^2 \lra \bbR$, allows us to use the notion of ends from the algebraic topology
(cf. \cite{HR}, see Definition~\ref{def-ends}) to analyze the behavior of level
sets at infinity and consequently to fully describe the topological structure of level
sets and their complements (which is achieved in Sections~\ref{sec-str} and
\ref{sec-nbr}).

Throughout the paper we use standard notation, as may be found in
\cite{BL,Kur1,Kur2,WhyAT}.

\begin{ack}
I would like to thank Professors W.B. Johnson, J. Lindenstrauss, R. Pol and G.
Schechtman for their interest in this work and many stimulating discussions,
O. Maleva for many constructive comments on a preliminary version of this paper, and
the organizers of the NSF Workshop in Linear Analysis and Probability at the Texas
A\&M University in College Station, Texas, where this work was started, for their
support and hospitality.
\end{ack}

\section{Number of components of level sets of uniform quotient
mappings from $\bbR^n$  to $\bbR$ is finite} \lb{sec-sep}

As a corollary of Theorem~\ref{finitecompl} using purely
topological arguments we will show that when $f:\bbR^n \lra \bbR$
is a uniform  quotient mapping then for each $t \in \bbR$, the
number of components of $f^{-1}(t)$ is finite
(Theorem~\ref{finite} below).  Our main tool is the following fact:

\begin{theorem} \lb{jan}
Let $B_{0}, B_{1} \subset {\mathcal{S}}_{n},\ n \geq 2$, be two
closed sets such that $B_0 \cap B_1 \subseteq \{q\}$ a singlepoint.  If
none of the sets $B_0$ or $B_1$ separates between points $p_1$ and
$p_2$ then their union $B_0 \cup B_1$ does it neither.
\end{theorem}

The above statement combines \cite [Theorem~59.II.11 and
61.I.7]{Kur2}  specialized to the situation in the present
paper.  In the case when $n=2$, Kuratowski refers to this fact as
the {\it first theorem of Janiszewski}, and it's general version is
called the {\it Phragmen-Brouwer theorem}.

Although the subject is closely related to some classical duality
theorems, cf. \cite{Newman,Spa85,BG99}, we were unable to find in the literature
results that we could directly use in the situation we deal with. We decided
to present a proof of the fact we needed, based on some standard arguments
concerning separation in $\mathbb{R}^n$.


We start
from two lemmas.

\begin{lemma}\lb{unique}
Let $A$ be an open connected subset of ${\mathcal{S}}_n$, so that
$\overline{A} \neq {\mathcal{S}}_n$  and $\Bd (A)= F_1 \cup F_2$
where $F_1, F_2$ are closed sets with $F_1 \cap F_2 \subseteq \{q\}$ a
singlepoint. Let $p_1 \in A$ and $p_2 \notin \overline{A}$.  Then
exactly one of the sets $F_1$ or $F_2$ separates between $p_1$ and $p_2$.
\end{lemma}

\begin{proof}
By Theorem~\ref{jan} we conclude that at least one of the sets
$F_1$ or $F_2$ separates $p_1$ and $p_2$. Suppose now that each of
$F_1$ and $F_2$ separates between $p_1$ and $p_2$.  Then there exist
components $C_1, C_2$ of ${\mathcal{S}}_n \setminus F_1$, ${\mathcal{S}}_n
\setminus F_2$ respectively so that
\begin{equation*}
\begin{split}
p_1 &\in C_1 \cap C_2, \text{ and  thus\ } A \subset C_1 \cap C_2,\\
p_2 &\notin C_1 \cup C_2.
\end{split}
\end{equation*}

Then
$$
\Bd(C_1 \cup C_2) \subset \Bd (C_1) \cup \Bd (C_2) \subset F_1
\cup F_2.
$$

Let $x \in F_1 \setminus \{q\}$.  Then for every neighborhood
$V_x$ of $x$ we have $V_x \cap A \neq \emptyset $, since $x \in
\Bd (A)$. Thus $V_x \cap C_2 \neq \emptyset$ and $x \in
\overline{C}_2$. Since $x \notin F_2$ we conclude that $x \in C_2$
and therefore $x \notin \Bd (C_1 \cup C_2)$.  Similarly, if $y \in
F_2 \setminus \{q\}$ then $y \notin \Bd (C_1 \cup C_2)$.  Thus
$\Bd (C_1 \cup C_2) \subset \{q\}$ which contradicts the fact that
$p_2 \notin C_1 \cup C_2$.
\end{proof}

\begin{lemma}\lb{separate}
Let $A$ be an open connected subset of ${\mathcal{S}}_n$ so that
$A \neq {\mathcal{S}}_n$ and $\Bd(A)= F_1 \cup F_2$ where $F_1,
F_2$ are closed sets with $F_1 \cap F_2 \subseteq \{q\}$ a singlepoint.
Suppose that ${\mathcal{S}}_n \setminus F_1$ is connected.  Then
for every $x \in F_1 \setminus \{q\}$ there exists a neighborhood
$U_x$ of ${x}$ so that $U_x \subset \overline{A}$.
\end{lemma}

\begin{proof}
Let $x \in F_1 \setminus \{q\}$.  Since $F_2$ is closed, there
exists a connected neighborhood $U_x$ of $x$ so that $U_x \cap F_2
= \emptyset$. Since $x \in \Bd (A)$, there exists $y \in U_x$ so
that $y \in A$.  Suppose that $U_x \setminus \overline{A} \neq
\emptyset$ and let $z \in U_x \setminus \overline{A}$.  Then $\Bd
(A)$ separates between the points $y$ and $z$.  But ${\mathcal{S}}_n
\setminus F_1$ is connected so $F_1$ does not separate between $y$ and
$z$.  Thus by Theorem~\ref{jan} we conclude that $F_2$ separates
between $y$ and $z$.  But this is a contradiction since $y,z$
belong to a connected set $U_x$ which is disjoint with $F_2$.
\end{proof}

With these tools we are ready to prove the main theorem of this section.

\begin{theorem} \lb{finite}
Let $f:\bbR^{n} \lra \bbR$ be a uniform quotient map.
  Then, for any $t\in \bbR$, a number of connected components of
$f^{-1}(t)$  is
finite and bounded by a function of $n$, $ \omega(\cdot)$ and $\Omega
(\cdot)$ only.
\end{theorem}

\begin{proof}
We consider $ \mathbb{R}^n$ as embedded in its one point compactification ${\mathcal{S}}_n$.
Denote
$K=f^{-1}(t)$.
By \cite[Lemma~5.2]{JLPS00}, $K$ is unbounded and therefore
 the closure of $K$ in ${\mathcal{S}}_n $ equals
$K\cup\{\infty\}$, and the closure  in ${\mathcal{S}}_n $ of every component of $K$ contains
$\{\infty\}$.
By Theorem~\ref{finitecompl}, $\bbR^n\setminus K$ and therefore also
${\mathcal{S}}_n \setminus \overline{K}$ has a
finite number of components, say
$$
{\mathcal{S}}_n \setminus \overline{K} = \bigcup_{j=1}^{m} C_j.
$$

Here $C_j\subset{\mathcal{S}}_n \setminus\{\infty\}$, so each $C_j$ can also be
considered as a subset of $ \mathbb{R}^n$. Note that $C_j$ cannot be bounded in $\mathbb{R}^n$,
so $\infty\in\Bd(C_j)\subset{\mathcal{S}}_n $ for all $j$.
  Suppose that there exists $j$, say $j=1$, so that $\Bd(C_1)$ has $m$ or more
connected components in $\mathbb{R}^n$. Then $\Bd(C_1)\subset \bbR^n$ can be presented as a
sum of $m $ disjoint closed sets in $\bbR^n$, which are not necessarily connected.
Thus after taking closures in ${\mathcal{S}}_n$ we see that
$$
\Bd(C_1) = F_1 \cup \cdots \cup F_m,
$$
where $\{F_k\}^m_{k=1}$ are closed sets in ${\mathcal{S}}_n$, not necessarily connected,
 so that $F_k \cap F_l \subseteq
\{\infty\}$ for all $k \neq l$.

Let $\{p_j\}^m_{j=1}$ be a collection of points such that $p_j \in
C_j$ for $j=1, \cdots, m$.
Since for each $j=2, \dots, m$, $\Bd (C_1)$ separates between $p_1$
and $p_j$, by Theorem~\ref{jan}, there exists $\sigma(j) \in
\{1,\cdots,m\}$ so that $F_{\sigma(j)}$ separates between $p_1$
and $p_j$.  By Lemma~\ref{unique}, $\Bd(C_1)\setminus F_{\sigma
(j)}$ does not separate between $p_1$ and $p_j$, so the choice of
$\sigma(j)$ is unique.  Thus  $\card(\{\sigma(j)\}^m_{j=2})\leq
m-1$. Hence there exists $j_0 \in \{1, \cdots, m\}$ so that
$F_{j_0}$ does not separate between $p_1$ and $p_i$ for all $i=2,
\cdots, m$.  Thus ${\mathcal{S}}_n \setminus F_{j_0}$ is
connected, and by Lemma~\ref{separate} for every $x \in F_{j_0}
\setminus \{\infty\}$ there exists a neighborhood $U_x$ of $x$ so
that $U_x \subset \overline{C}_1$.  But $f(C_1)\subset (t,
\infty)$ or $f(C_1) \subset(-\infty, t)$, thus $f(U_x) \subset [t,
\infty)$ or $f(U_x) \subset (-\infty, t]$, which contradicts the
fact that $f(U_x) \supset B(f(x),\e) = (t-\e,t+\e)$ for some $\e > 0$.
This contradiction yields that $\Bd (C_1)$ has at most $(m-1)$
components in $\bbR^n$. Similarly, for every $j\in\{1,\dots,m\}$,
$\Bd(C_j)$ has at most $(m-1)$ components and
 since every component of $K$ contains a
component of $\Bd(C_j)$ for at least one $j \in \{1,\cdots, m\}$,
we conclude that the number of components of $K$ is smaller or equal
than $m(m-1)$.
\end{proof}

\begin{cor} \lb{sepplane}
Let $f: \bbR^n \lra \bbR$ be a uniform quotient mapping.  Then, for
each   $t \in \bbR$, each component of $f^{-1}(t)$ separates
$\bbR^n$.
\end{cor}

This Corollary has word for word the same proof as
\cite[Proposition~5.4]{JLPS00},  since by Theorem~\ref{finite}, for
each $t \in \bbR$, $f^{-1}(t)$ has a finite number of components.

\section{Local connectedness of level sets}
\lb{sec-lc} In this section we show that all level sets of
co-Lipschitz uniformly continuous mappings from $\bbR^2$   to $\bbR$ are
hereditarily locally connected. This is a very strong property
which will enable us to give a  detailed description of the
structure of the level sets, see Sections~\ref{sec-str} and
\ref{sec-nbr}.

As mentioned in the Introduction, we do not know whether there exist co-Lipschitz uniformly
continuous maps or Lipschitz quotient maps from $\bbR^3$ to $\bbR$, or
in general from $\bbR^n$ to $\bbR^k$, with non-locally connected level sets.
However we do know that there exist uniform quotient maps from  $  \bbR^2$ to  $\bbR$ which have non-locally
connected level sets (see Section~\ref{example}).

We begin by recalling some basic definitions.

\begin{defn}
A topological space  $S$ is
said to be {\it locally connected at a point x} if for every open set $U$ containing
$x$ there is a connected open set $V$ so that $x \in V \subset U$.  The space $S$ is
{\it locally connected} if it is locally connected at each point and $S$ is {\it
hereditarily locally connected} if every subcontinuum of $S$ is locally connected.
\end{defn}

We
will use the following characterization of hereditary local connectedness:

\begin{theorem} \lb{hlc} \cite[V.(2.1) and I.(12.2)]{WhyAT}
A locally compact connected set $S$ is hereditarily locally connected if and only  if
$S$ does not contain a continuum of convergence.
\end{theorem}

Recall   that if a continuum $K$ is a subset of a set $M$ then $K$
is called a {\it continuum of convergence} of $M$ provided that there exists in $M$ a
sequence of mutually exclusive continua $K_1, K_2, \dots,$ no one of which contains a
point of $K$ and which converges to $K$ as a limit, i.e. $K \cap
\bigcup_{i=1}^{\infty} K_i = \emptyset$ and $\lim[K_i]_i = K$.

Here $\lim[K_i]_i$ denotes the limit of a sequence $[K_i]_i$ which is defined
as follows (cf. \cite[Section I.7]{WhyAT} or \cite[Chapter 11, Section 29]{Kur1}): The
set of all points $x$ such that every neighborhood of $x$ contains points of
infinitely many sets of $[K_i]_i$ is called the {\it limit superior} of $[K_i]_i$ and
is denoted $\lim \sup [K_i]_i$.  The set of all points $y$ such that every
neighborhood of $y$ contains points of all but a finite number of the sets $[K_i]_i$
is called the {\it limit inferior} of $[K_i]_i$ and is denoted $\lim \inf [K_i]_i$.  If
$\lim \sup [K_i]_i = \lim \inf [K_i]_i$ then we say that the collection $[K_i]_i$ is
{\it convergent} and we write $\lim[K_i]_i = \lim\sup [K_i]_i = \lim \inf [K_i]_i$ and
we call $\lim[K_i]_i$, the {\it limit} of $[K_i]_i$.

We will prove that if $f: \bbR^2 \lra \bbR$ is a co-Lipschitz uniformly continuous mapping then for
all $t \in \bbR, f^{-1}(t)$ does not contain a continuum of convergence.  For this we
will need the following ``bottleneck lemma'', whose proof is very similar to the proof
of \cite[Lemma~5.3]{JLPS00}.

\begin{lemma}\lb{bottleneck}
Let $f:\bbR^2 \lra \bbR$ be a co-Lipschitz uniformly continuous map with co-Lipschitz constant
 $1$ and a modulus of uniform continuity $\Om$.
Let $K_1, K_2$ be disjoint subcontinua of $f^{-1}(0)$ and let $\alpha \in \bbR_+$.  If
there exist points $x_1, x_2 \in K_1$, $y_1,y_2 \in K_2$ so that, for $i=1,2$,
$$
d (x_i, y_i)\leq \alpha,
$$
then
$$d(x_1, x_2)\leq 2\Om(\frac\al2)+4\alpha.$$
\end{lemma}

For the proof we will need   the following  basic lemma concerning the lifting of
Lipschitz curves which was established in \cite{BJLPS}.

\begin{lemma} \lb{lift} \cite[Lemma~4.4]{BJLPS} Suppose that $f:\bbR^{n}\lra X$ is
 continuous and
co-Lipschitz with constant one, $f(x)=y$, and $\xi:[0,\infty) \lra X$
is a curve with Lipschitz constant one, and $\xi(0)=y$.  Then there
is a curve $\phi:[0,\infty)\lra \bbR^{n}$ with Lipschitz constant one
such that $\phi(0)=x$ and $f(\phi(t))=\xi(t)$ for $t\ge 0$.
\end{lemma}

\begin{proof}[Proof of Lemma~\ref{bottleneck}]
If $d(x_1,x_2)\leq 2 \alpha$ then we are done, so assume without loss of generality
that $d(x_1, x_2)> 2\alpha$. For $i=1,2$, let $I_i$ be the segment connecting $x_i$ and
$y_i$, i.e.
$
I_i = \{(1-t)x_i + ty_i \ : \ t \in[0,1]\}.
$
Then  $length(I_i)\leq \alpha$, for $i=1,2$, and thus, if $d(x_1,x_2)>2\alpha$ then
$I_1 \cap I_2 = \emptyset$. Set
$$t_i = \sup \{t \in [0,1] : (1-t)x_i + ty_i \in K_1\}.$$
Since $K_1$ and $I_i$ are
compact and $y_i \notin K_1$ we get that $t_i \in [0,1)$.  Define
$$
\overline{x_i} \overset{ \text{def}}{=} (1-t_i)x_i + t_i y_i \in K_1.
$$
Now set
$$
s_i = \inf \{t \in [t_i,1]: (1-t)x_i + ty_i \in K_2\}.
$$
Similarly as above, since $K_2$ is compact and $\overline{x_i} \notin K_2$, we get
that $s_i \in (t_i,1]$.  Define
$$
\overline{y_i} \overset{ \text{def}}{=} (1-s_i)x_i + s_i y_i \in K_2.
$$

Further, for $i=1,2,$   define segments with endpoints $\overline{x_i}, \overline{y_i}$,
$$
J_i \overset{ \text{def}}{=} \{(1-t)\overline{x_i} + t\overline{y_i} : t \in [0,1]\}.
$$

Then we get that $J_i \cap K_1 = \{ \overline{x_i}\}$, $J_i \cap K_2 = \{
\overline{y_i}\}$  and $J_i \cap J_2 = \emptyset$ (since $J_i \subset I_i$ which were
disjoint).
Further
\begin{equation}\lb{alpha}
d(\overline{x_i}, \overline{y_i}) \leq d(x_i, y_i)< \alpha.
\end{equation}

By \cite[Theorem~62.V.6]{Kur2} there exists an open connected region $G$ whose
boundary is contained in $K_1 \cup K_2 \cup J_1 \cup J_2$.
Since $K_1 \cup K_2 \subset f^{-1} (0)$, and by \eqref{alpha}, we conclude that for
all $x \in \Bd (G)$,
\begin{equation}\lb{bdry}
|f(x)| \leq \Om(\frac\al2).
\end{equation}

Let $x_0 \in G$ be such that for $i=1,2$
$$
d(x_0, \overline{x_i}) \geq \frac{1}{2}d (\overline{x_1},\overline{x_2}).
$$

Such a point $x_0$ exists in $G$ since $G$ is open and connected and thus $G$ is
path-connected.  By Lemma~\ref{lift} there exists a curve $\phi:[0, \infty)\lra \bbR^2$
with Lipschitz constant one, $\phi (0)= x_0$ and $f(\phi (t))= f (x_0)+ t \ sign
(f(x_0)).$  Since this curve is clearly unbounded, there exists $\tau >0$ so that
$\phi (\tau) \in \Bd (G)$. Then, by \eqref{bdry} and since $\phi$ is Lipschitz with
constant one,
\begin{equation*}
\begin{split}
\Om(\frac\al2) &\geq |f (\phi(\tau))| \geq \tau \geq \| \phi
(\tau)- \phi (0) \|= \| \phi (\tau)- x_0 \| \\
& \geq d (x_0, J_1 \cup  J_2)\geq \min_{i=1,2} (d(x_0, \overline{x_i}))-\alpha
\\
&\geq
\frac{1}{2}d(\overline{x_1},\overline{x_2})-\alpha.
\end{split}
\end{equation*}
Thus
$$
d(\overline{x_1},\overline{x_2}) \leq 2\Om(\frac\al2)+2\al,
$$
and
$$
d(x_1, x_2)\leq d(x_1,\overline{x_1})+ d(\overline{x_1},\overline{x_2})
+ d(\overline{x_2},x_2) \leq 2\Om(\frac\al2)+4\al.
$$
\end{proof}

\begin{prop}\lb{r2}
Let $f:\bbR^2 \lra \bbR$ be a co-Lipschitz uniformly continuous map.  Then for every  $ t \in \bbR,
f^{-1}(t)$ is hereditarily locally connected.
\end{prop}

\begin{proof}
Without loss of generality we will assume that $f$ is a co-Lipschitz uniformly continuous map
with a co-Lipschitz
constant  $1$, and that $t=0$.  By Theorem~\ref{hlc}, it is enough to
show that $f^{-1}(0)$ does not contain a continuum of convergence.

Suppose for contradiction that $K_0$ is a continuum of convergence in $f^{-1}(0)$  and
let $x_1, x_2 \in K_0$, and $ \beta \overset{ \text{def}}{=} d(x_1, x_2)>0$.  Let
$[K_i]^\infty_{i=1}$ be the sequence of mutually disjoint subcontinua of $f^{-1}(0)$
with $\bigcup_i K_i \cap K_0 = \emptyset$ and $\lim[K_i]_i = K_0$.  Then, by the definition
of the limit (see also \cite[Theorem I.(7.2)]{WhyAT}), for every $\e > 0$ there
exists $n \in \bbN$ so that for all $x \in K_0$
$$
d(x, K_n)< \e.
$$
Thus for $i=1,2$ there exists $y_i \in K_n$ with
$$
d(x_i, y_i)< \e.
$$
Hence, by Lemma~\ref{bottleneck}, $d(x_1,x_2) \leq  2\Om(\e/2)+4\e $. Since $\e$
is   arbitrary and $\lim_{r\to 0}\Om(r)=0$,
  we conclude that $d(x_1,x_2)=0$ which    contradicts
the fact that $K_0$ is a nontrivial subcontinuum.
\end{proof}

\section{First description of the structure of level sets}\lb{sec-str}

The aim of this section is to   obtain  a characterization
of the form of any
  closed, hereditarily locally connected, locally compact, connected set with no end
 points and containing no simple closed curve, and to apply it to
 describe the structure of level sets of co-Lipschitz uniformly continuous mappings $f: {\mathbb{R}}^2
\longrightarrow  {\mathbb{R}}$ (Theorem~\ref{graph} and Remark~\ref{rem-d}).  For that we will
need the notions of a dendrite, an order of a point, an end
point and a cut point of a topological space $M$.  We recall
their definitions below.

\begin{defn} \cite[Chapter III]{WhyAT}, \cite[Chapter~VI, \S51]{Kur2}\lb{deford}
Let $M$ be a space and $\mathfrak{n}$ a cardinal number.  We say that a point
$x \in M$ is of {\it order $\leq \mathfrak{n}$} in $M$ provided that for any
neighborhood $V$ of $x$ in $M$, there exists a neighborhood $U$ of
$x$ in $M$ with $U \subset V$ and $\card(\Bd (U)) \leq \mathfrak{n}$.

A point $x \in M$ which is of order one in $M$ will be called an
{\it end point} of $M$.
\end{defn}

\begin{defn} \cite[Section III.1]{WhyAT},
\cite[Definition~47.VIII.2]{Kur2}
If $M$ is a connected set and
$p$ is a point of $M$ such that the set $M\setminus\{p\}$ is not
connected, then $p$ is called a {\it cut point} of $M$.
\end{defn}

\begin{defn}  \cite[Section ~V.1]{WhyAT} A
continuum $M$ is called a {\it dendrite} (or an {\it acyclic curve})
provided that $M$ is locally connected and contains no simple
closed curve.
\end{defn}

Dendrites constitute a very important class of
continua, and they have been extensively studied. We recall here
a couple of important properties of dendrites, that we will use.

\begin{theorem} \lb{cd} \cite [V.(1.1) and V.(1.2)] {WhyAT}
Let $M$ be a continuum.  The following statements are equivalent:

(1) $M$ is a dendrite.

(2) Every point of $M$ is either a cut point or an end point.

(3) $M$ is locally connected and one and only one arc exists between
any two points in $M$.
\end{theorem}

We observe here a simple property of end points which we state in
the form of a lemma for an easy reference:

\begin{lemma} \lb{sep-end}
Let $M$ be a closed connected, locally connected subset of
$\bbR^{n}$.  Suppose that $p$ is an end point of the subset $B
\subset M$ such that:
\begin{itemize}
\item[(a)] there exists an open set $U \subset \bbR^{n}$ so that $p \in
U$ and $M \setminus B \subset \bbR^{n} \setminus \overline{U}$; or
\item[(b)] $B$ is a component of $M \setminus A$ for some subcontinuum
$A$ of $M$.
\end{itemize}

Then $p$ is an end point of $M$.
\end{lemma}

\begin{proof}
For (a) let $V$ be any neighborhood of $p$.  Then, since $p$ is an
end point of $B$, there exists an open set $V_{1} \subset V \cap
U$ so that $p \in V_{1}$ and $\card(\Bd(V_{1}) \cap B)=1$.  Since
$M \setminus B \subset \bbR^{n} \setminus \overline{U}$ we get
that $(M \setminus B) \cap \overline{V_{1}}=\emptyset$ and thus
$\Bd(V_{1}) \cap M=\Bd(V_{1}) \cap B$, which ends the proof of part (a).

For (b) let $V$ be any neighborhood of $p$.  Since $M$ is locally
connected, there exists an open set $V_{1} \subset V \setminus A$
so that $V_{1} \cap M$ is connected.  Since $p \in B$ and $B$ is a
connected component of $M \setminus A$ we conclude that $V_{1}
\cap M = V_{1} \cap B$.  Thus, since $p$ is an end point of $B$,
there exists an open set $V_{2} \subset V_{1}$ so that $p \in
V_{2},\ \overline{V_{2}} \subset V_{1},\ \card(\Bd(V_{2}) \cap B)=1$
and $\Bd(V_{2})\cap B = \Bd(V_{2}) \cap M$, which ends the proof of
part (b).
\end{proof}

Our first observation concerning the structure of level sets
of co-Lipschitz uniformly continuous mappings
is the following:

\begin{cor} \lb{dendrite}
Let $f : \bbR^{2} \lra \bbR$ be a co-Lipschitz uniformly continuous map.
Then for every $t \in \bbR$ and every subcontinuum $M$ of
$f^{-1}(t)$, $M$ is a dendrite.
\end{cor}

\begin{proof}
It is easy to see that when $f$ is a co-Lipschitz uniformly continuous map then
for all $t$, $f^{-1}(t)$ cannot contain a simple closed curve.
Indeed, since $f$ is   co-Lipschitz, $f^{-1}(t)$ has
empty interior and if a simple closed curve $C$ was contained in
$f^{-1}(t)$ then $\bbR^{2} \setminus f^{-1}(t)$ would have a
bounded component $A$ contained in the region inside the curve
$C$.  But then $f(A)=(t,\infty)$ or $(-\infty,t)$, which is
impossible since $\overline{A}$ is compact and $f$ is continuous.

Further, by Proposition~\ref{r2}, every subcontinuum $M$ of
$f^{-1}(t)$ is locally connected.  Thus $M$ is a dendrite.
\end{proof}

Our next goal is to show that every $f^{-1}(t)$ is of a
particularly simple form, that every point of $f^{-1}(t)$ is of
finite order and only finitely many points in $f^{-1}(t)$ have
order bigger than $2$.  Thus we will show that $f^{-1}(t)$ has a
graph structure.  Recall the following

\begin{defn} (cf. e.g. \cite[Section~X.1]{WhyAT})
A set $A$ is called a (finite linear) {\it graph} provided $A$ is
the union of a finite set $V$ of points, called {\it vertices},
and a finite number of open mutually disjoint arcs
$\alpha_{1},\alpha_{2},\dots,\alpha_{n}$, called {\it edges}, so
that the two end points of each edge $\alpha_{i}$ are distinct and
belong to $V$.  A graph which contains no simple closed curve is
called a {\it tree} or an {\it acyclic graph} (see e.g.
\cite[Definition~9.25]{Nadler}).
\end{defn}

We start from the following:

\begin{prop} \lb{endpoint}
Let $f:\bbR^{2}\lra\bbR$ be a co-Lipschitz uniformly continuous map and
let $K$ be a component of $f^{-1}(t)$ for some $t \in \bbR$.  Then
every point of $K$ is a cut point, i.e. $f^{-1}(t)$ has no end
points.
\end{prop}

\begin{proof}
Suppose that $x \in K$ is not a cut point of $K$.  Then by
Corollary~\ref{dendrite} and Theorem~\ref{cd}, $x$ is an endpoint
of $K$.  It follows from \cite[Theorem~26]{Why27} (cf. also
\cite[Proof of Theorem~27]{G28}, where this fact is attributed to
R.G. Lubben), that if $x \in K$ is an endpoint of $K$ then $x$
belongs to the boundary of exactly one component of $\bbR^{2}
\setminus K$.  But then, since $\bbR^{2} \setminus f^{-1}(t)$ and
thus also $\bbR^{2} \setminus K$ have finite number of components,
there exists a neighborhood $U$ of $x$ so that $U$ intersects
exactly one component of $\bbR^{2} \setminus K$.  Hence $f (U)
\subset (t,\infty)$ or $f(U) \subset (-\infty,t)$, which
contradicts the fact that $f$ is co-Lipschitz.
\end{proof}

To finish the analysis of the structure of level sets $f^{-1}(t)$
we will need one more notion -- the notion of a number of ends of an
unbounded locally connected set.

\begin{defn}\lb{def-ends} \cite[Definition~1.18]{HR}
We say that a connected locally connected Hausdorff space $W$ {\it has  at least $k$
ends}  if there exists an open subspace $V \subseteq W$ with compact closure
$\overline{V}$ so that $W \setminus \overline {V}$ has at least $k$ unbounded
components.  The space $W$ {\it has exactly k ends} if $W$ has at least $k$ ends but
not at least $k+1$ ends.  If $W$ has exactly $k$ ends we will write $\# e(W)=k$.
\end{defn}

One should be careful not to confuse ends with end points.  We
think of ends, intuitively, as infinite ends of unbounded sets.
In fact, there exist ways of making this intuition precise, by defining ends using homotopy
classes of unbounded paths contained in the space $W$ (see \cite{HR}), but we will not need
this for our present purpose.

Clearly continua never have any ends, but unbounded locally
connected sets may have some end points in addition to the fact
that they always have at least one end.

If a locally connected space $W$ has a finite number of connected components,  $W =
\bigcup^m_{j=1}C_j$, then we will use notation $\#e(W)$ to mean the sum of $\#e(C_j)$, i.e.
$$
\#e(W) \overset{ \text{def}}{=} \sum^m_{j=1}\#e (C_j).
$$

We note here that it follows from the local connectivity of $W$, by
\cite[Theorem~3-9]{HY},  that if $V$ is an open subset of $W$ with compact closure,
then $W\setminus \overline{V}$ has at most a finite number of unbounded
components.  If $W$ has exactly $k$ ends then there exists an open subspace $V
\subseteq W$ with compact closure so that $W \setminus \overline {V}$ has exactly
$k$ unbounded components.  Moreover we have the following:

\begin{prop}\lb{endstr} \cite[Proposition~1.20 and its proof]{HR}
For an unbounded connected locally connected closed space $W\subset\bbR^n$ with
 exactly $k$ ends there exists an open set $U\subset\bbR^n$ with compact closure so that
$W$ can be expressed as
$$
W = W_0 \cup \bigcup_{j = 1}^k \overline{W(j)},
$$
where $W_0 = W\cap \overline{U}$ is   connected and compact, sets $W(j)$, for $j=1,\dots,k,$
 are
   connected components of
$ W\setminus U$ and each $W(j)$ has  exactly one end.
\end{prop}

As a corollary we obtain the main theorem of this section which
describes the structure of level sets $f^{-1}(t)$.

\begin{theorem} \lb{graph}
Let $f:\bbR^{2}\lra\bbR$ be a co-Lipschitz uniformly continuous map.  Then
for every $t \in \bbR$, every component $K$ of $f^{-1}(t)$ has a
representation of the form:
$$
K=K_{0} \cup \bigcup^{n}_{j=1} K_{j},
$$
where $n \in \bbN,\ n =\#e(K),\ K_{0}$ is a compact connected tree with exactly $n$ endpoints,
 each $K_{j}$ is a ray, that is a closed unbounded set
homeomorphic with $[0,\infty)$, sets $\{K_{j}\}^{n}_{j=1}$ are
mutually disjoint, for all $j$, $\card(K_{j} \cap K_{0})=1$   and the unique point
in the intersection $K_{j} \cap K_{0}$ is an end point of $ K_{0}$ and of $K_{j}$.
\end{theorem}

\begin{defn} We will use the term {\it unbounded finite graph} for the sets of the form
described
in Theorem~\ref{graph}.
\end{defn}

\begin{proof}
Let $K$ be a component of $f^{-1}(t)$ for some $t \in \bbR$.  By
Proposition~\ref{r2} $K$ is locally connected and thus, by
\cite[Theorem~3-9]{HY}, $K$ has exactly $n$ ends for some $n \in
\bbN ,\ n \geq 1$.  Thus, by Proposition~\ref{endstr}, there exists
an open set $U \subset \bbR^{2}$ with compact closure so that $K$
can be expressed as
$$
K=W_{0} \cup \bigcup^{n}_{j=1} \overline{W(j)},
$$
where $W_{0}=K \cap \overline{U}$ is connected and compact, and
sets $W(j) \subset K \setminus U$ are connected components of $K
\setminus U$ and each $W(j)$ has exactly one end.
We define for   $j \in
\{1,\dots,n\}$,
$$K_j\DEF\overline{W(j)}.$$

First we notice that, since every subcontinuum of $K$ is a
dendrite, we have for all $j=1,\dots,n$:
\begin{equation} \lb{unique1}
\card(\overline{W(j)} \cap W_{0})=1.
\end{equation}

To see this, assume, for contradiction, that for some $j \in
\{1,\dots,n\}$, there exist two points $x_{1} \neq x_{2}$ in
$\overline{W(j)} \cap W_{0}$.  Let $V_{1},V_{2}$ be open disjoint
neighborhoods of $x_{1},x_{2}$, respectively, so that sets $V_{i}
\cap \overline{W(j)}$, for $i=1,2$, are arcwise connected; this is
possible since $\overline{W(j)}$ are locally connected.  For
$i=1,2$, let $y_{i} \in V_{i} \cap W(j)$ and let
$\eta_{i}:[0,1]\lra V_{i}\cap\overline{W(j)}$ be a path so
that $\eta_{i}(0)=y_{i},\ \eta_{i}(1)=x_{i}$.  Set, for $i=1,2$,
$$
\begin{aligned}
t_{i} &\DEF \inf \{t \in [0,1]:\eta_{i}(t)\in W_{0}\},\\
z_{i} &\DEF \eta_{i}(t_{i}) \in \overline{W(j)} \cap W_{0},\\
\widehat{\eta_{i}} &\DEF \eta_{i}([0,t_{i}]).
\end{aligned}
$$

Next, note that since $W(j)$ is open and locally connected, by
\cite[Theorem~50.I.2]{Kur2}, there exists an arc $\alpha \subset
W(j)$ with endpoints $y_{1}$ and $y_{2}$.  Now, let $\beta =
\widehat{\eta_{1}} \cup \alpha \cup \widehat{\eta_{2}}$.  Then $\beta$
in an arc with endpoints $z_{1}$ and $z_{2}$, and $\beta \setminus
\{z_{1},z_{2}\} \subset W(j)$.

On the other hand, since $z_{1},z_{2} \in W_{0}$ and since, by
Corollary~\ref{dendrite}, $W_{0}$ is a dendrite we conclude, by
Theorem~\ref{cd}(3), that there exists an arc $\gamma \subset
W_{0}$ with endpoints $z_{1}$ and $z_{2}$.

Thus we obtained two arcs $\beta$ and $\gamma$ in $K$ with $\beta
\cap \gamma = \{z_{1},z_{2}\}$ which contradicts
Theorem~\ref{cd}(3) and ends the proof of \eqref{unique1}.

 Define, for $j=1,\dots,n$, $w_{j}$ to be the unique point, by \eqref{unique1}, of the
intersection
$\overline{W(j)} \cap W_{0}$.
We claim that
\begin{equation}\lb{unique2}
\overline{W(j)}=W(j) \cup \{w_{j}\}.
\end{equation}

Indeed, since $W(j)$ is a component of $K \setminus U,\ \Bd(W(j))
\subset \Bd(U)$ that is:
$$
\overline{W(j)} \subset W(j) \cup \overline{U}.
$$
Since $K$ is closed and $W_{0}=K \cap \overline{U}$ we get
$$\overline{W(j)} \subset K \cap (W(j) \cup \overline{U}) = W(j)
\cup W_{0}.$$
Thus \eqref{unique2} follows   from \eqref{unique1} and
the definition of $w_{j}$.

Now fix a point $w_{0} \in W_{0} \cap U$.
By Corollary~\ref{dendrite} and Theorem~\ref{cd}(3), for each $j
\in \{1,\dots,n\}$, there exists exactly one arc $\sigma_{j} \subset
W_{0} \subset K$ with endpoints $w_{0}$ and $w_{j}$.  Set
$$
K_{0} \DEF \bigcup^{n}_{j=1} \sigma_{j}.
$$

It is clear from the definition that $K_{0}$ is a compact
connected graph with exactly $n$ endpoints $ \{w_{1},\dots,w_{n}\}$,
and, by Corollary~\ref{dendrite}, $K_{0}$ does not
contain any simple closed curve, thus $K_{0}$ is a tree.
Also, by \eqref{unique2}, for all $j
\in \{1,\dots,n\}$, $\overline{W(j)} \cap K_{0}=\{w_j\}$, an end point of $K_0$.

We will now show that $W_{0}=K_{0}$.

Clearly $K_{0} \subset
W_{0}$, so suppose, for contradiction, that $W_{0} \setminus
K_{0}$ is nonempty and let $A$ be the closure of a connected
component of $W_{0} \setminus K_{0}$.
By \cite [Theorem~6.8]{Nadler} there exists a non-cut point $p$ of $K_{0}
\cup A$ so that $p \in A$ and $p$ is a non-cut point of $W_{0}$,
thus $p$ is an end point of $W_{0}$.  Since $p \in A$, we see that
$p \notin K_{0}$ and $p \notin \{w_{1},w_{2},\dots,w_{n}\} = K_{0}
\cap \bigcup^{n}_{j=1} \overline{W(j)}$.  By \eqref{unique2}
$$
W_{0} \cap \bigcup^{n}_{j=1} \overline{W(j)} \subseteq
\{w_{1},w_{2},\dots,w_{n}\}.
$$

Thus
$$
p \notin \bigcup^{n}_{j=1} \overline{W(j)},
\text{\ \ \ \ \
and
\ \ \  }
p \in U.
$$

Hence, by Lemma~\ref{sep-end}(a), since $p$ is an end point of
$W_{0}$, $p$ is also an end point of $K$, which is a contradiction
with Proposition~\ref{endpoint}.  Thus we have shown that
$W_{0}=K_{0}$.

To finish the proof of the theorem it is only left to show that
for each $j \in \{1,\dots,n\}$, the set $W(j)$ is homeomorphic with
the real line.  To see this we will use the classical theorem of
Ward \cite{W36}, which characterizes the real line as a non-empty
connected, locally connected, separable metric space which is cut
by each of its points into exactly $2$ components.

Clearly, each $W(j)$ is a non-empty, connected, locally connected,
separable metric space; and, by Proposition~\ref{endpoint}, every
point of $W(j)$ is a cut point of $K$, and thus also of $W(j)$.
Suppose that there exists a point $x \in W(j)$ so that $W(j)
\setminus \{x\}$ has $3$ or more components.  Since $W(j)$ has
exactly one end, exactly one component of $W(j) \setminus \{x\}$
is unbounded.  Thus $W(j) \setminus \{x\}$, and therefore also
$\overline{W(j)} \setminus \{x\} = (W(j) \cup \{w_{j}\}) \setminus
\{x\}$, have at least $2$ bounded components, say $C_{1},C_{2}$.

Let $V$ be an open set with compact closure so that
$\overline{C_{1}} \cup \overline{C_{2}} \subset V$.  Then
$C_{1},C_{2}$ are components of $\overline{W(j)} \cap
\overline{V}$.  Without loss of generality we assume that $w_{j}
\notin C_{1}$.  Then, by \cite [Corollary~5.9]{Nadler},
$$
\overline{C_{1}}=C_{1} \cup \{x\},
$$
and, by \cite[Theorem~6.8]{Nadler}, there exists a non-cut point $p$ of
$\overline{C_{1}}$ so that $p \not= x$ and $p$ is a non-cut point
of $\overline{W(j)} \cap \overline{V}$.  Hence, by
Corollary~\ref{dendrite} and Lemma~\ref{sep-end}(a), $p$ is an end point
of $\overline{W(j)}=W(j) \cup \{w_{j}\}$.  Since $p \not=
 w_{j}$, by Lemma~\ref{sep-end}(b), $p$ is an end point
of $K$ and
we get a contradiction with Proposition~\ref{endpoint}.  Thus we
have shown that every point of $W(j)$ cuts $W(j)$ into exactly $2$
components and hence $W(j)$ is homeomorphic to the real line and
$\overline{W(j)}=W(j) \cup \{w_{j}\}$ is a ray, which ends the
proof of the theorem.
\end{proof}

\begin{rem} \lb{rem-d}
As the reader has surely noticed,   the above proof and hence also
the conclusion of Theorem~\ref{graph} is valid for any set $K$ such that:
\begin{quotation}
{\it $K\subset\bbR^m$, $m\ge 2$,
 is a closed, hereditarily locally connected, locally compact, connected set with no end
 points and containing no simple closed curve.}
\end{quotation}
\end{rem}

\begin{rem} \lb{rem-tree}
There is an alternative way to prove Theorem~\ref{graph}.  Instead
of the fairly direct and self-contained proof presented above, one
could use results of Shimrat \cite{Sh55,Sh59} and Stone
\cite{S68} who (among others) studied the structure of sets whose
every point is a cut point.  In particular Shimrat \cite{Sh59}
fully described locally connected sets whose every point is a cut
point and this characterization when refined with \cite[Theorem~3-9]{HY}
and the assumption that the set is closed also gives the
statement of Theorem~\ref{graph}.

On the other hand, Stone \cite{S68} gave a characterization of
finite linear graphs, from which it follows easily that a dendrite
is a (finite linear) tree if an only if it has a fine number of
end points.  This result of Stone has been, using different
methods, reproved and strenghtened by Nadler \cite{N93} (cf. also
\cite[Theorem~9.24]{Nadler}).  It is
clear that this characterization is closely related with
Theorem~\ref{graph} and indeed it is possible to obtain a proof of
Theorem~\ref{graph} using these results.

However we felt that
following either of these two routes of reasoning would be
technically more complicated than the presented direct proof.
\end{rem}

\section {Number and form of components of level sets}
\lb{sec-nbr}

In this section we present an exact characterization of the form
of level sets of a co-Lipschitz uniformly continuous map $f$ from $\bbR^{2}$ to
$\bbR$  (Theorem~\ref{main}), which significantly refines
Theorem~\ref{graph}.  In particular, we obtain an affirmative
answer to the question posed in \cite{JLPS00} whether the number
of components of level sets $f^{-1}(t)$ or of $\bbR^{2} \setminus
f^{-1}(t)$ are constant after excluding finitely many values of
$t$.  We begin with the statement of our main characterization
theorem.

We will use the notation $\#c(W)$ to denote the number
of components of the set $W$.

\begin{theorem} \lb{main}
For every pair of a constant $c>0$ and a function $\Om(\cdot)$ with $\lim_{r\to 0}\Om(r)=0$,
there exists a natural number $M=M(c,\Om)$, so
that for every co-Lipschitz uniformly continuous map $f:\bbR^{2} \lra \bbR$ with a co-Lipschitz
constant $c$
and a modulus of uniform continuity $\Om$,
there exists a natural number $n=n(f) \leq M$ and a finite subset $T_{f}$ of
$\bbR$, with $\card(T_{f}) \leq n-1$, so that:
\begin{itemize}
\item[(1)]
for all $t \in \bbR$,
$$
\#e(f^{-1}(t))=2n,
$$
that is there exists $R_0\in \bbR$ so that for every $R>R_0$, $f^{-1}(t)\setminus B(0,R)$
has exactly $2n$ unbounded  components.

Moreover, if $\{C_{R,i}\}_{i=1}^{2n}$ are the   unbounded components of
$f^{-1}(t)\setminus B(0,R)$, then
for all $i\ne j$,
$$\lim_{R\to \infty} d(C_{R,i}, C_{R,j})=\infty;$$

\item[(2)]
for all $t \in \bbR \setminus T_{f}$,

\begin{itemize}
\item[(a)]
$\#c(f^{-1}(t))=n,$

\item[(b)]
$\#c(\bbR^{2}\setminus f^{-1}(t))=n+1,$

\item[(c)] each component of $f^{-1}(t)$ is
homeomorphic with the real line and separates the plane into
exactly 2 components;
\end{itemize}

\vspace{2mm}

\item[(3)] for all $t_{i} \in
T_{f}$,

\begin{itemize}
\item[(a)]
$\#c(f^{-1}(t_{i}))<n,$

\item[(b)]
$\#c(\bbR^{2}\setminus f^{-1}(t_{i}))=2n+1-\#c(f^{-1}(t_i)) \in
(n+1,2n),$

\item[(c)]
each component of $f^{-1}(t_{i})$ is an unbounded finite graph,
i.e. has the form described in Theorem~\ref{graph}
\end{itemize}
\end{itemize}
\end{theorem}

\begin{rem} \lb{rem-a}
Theorem~\ref{main}  is analogous to a result of Johnson, Lindenstrauss, Preiss and Schechtman
\cite{JLPS00}, who proved that for every  pair of a constant $c>0$ and a function $\Om(\cdot)$
with $\lim_{r\to 0}\Om(r)=0$,
there exists a natural number $M=M(c,\Om)$, so
that for every co-Lipschitz uniformly continuous map $f:\bbR^{2} \lra \bbR^2$ with a co-Lipschitz
constant $c$
and a modulus of uniform continuity $\Om$, there exists a natural number $n=n(f) \leq M$ and a polynomial $P_f$ with degree
  equal to $n$, so that $f=P_f\circ h_f$, where $h_f$ is a homeomorphism of
$\bbR^2$. Hence there exists a
finite set $T_{f}\subset\bbR^2$   with $\card(T_{f}) \leq n\leq M$, so that for all
$t\in\bbR^2\setminus T_f$, $\card(f^{-1}(t))=n$ and for all
$t_i\in  T_f$, $\card(f^{-1}(t_i))<n$, analogously with parts (2a) and (3a) of
Theorem~\ref{main}.

For Lipschitz quotient maps from $\bbR^{2}$ to $\bbR^2$,
Maleva \cite{M02} studied the dependence of the number $M(c,L)$ on the Lipschitz and
co-Lipschitz constants $L$ and $c$. Maleva proved in particular that there exists
 a scale $0<\dots<\ro_2^{(m)}<\dots<\ro_2^{(1)}<1$ such that for any
 Lipschitz quotient mapping $f:\bbR^{2} \lra \bbR^2$ the condition $c/L>\ro_2^{(m)}$
implies that $\card(f^{-1}(t)) \leq m$  for all
$t\in\bbR^2$ (in fact this holds with $\ro_2^{(m)}=1/(m+1)$) \cite[Theorem~2]{M02}.

It is natural to ask whether a similar scale exists for the numbers $M(c,\Om)$ defined
in Theorem~\ref{main}.
[After reading a preliminary version of this paper,  Maleva proved the existence of such
a scale. Namely she proved that if $c/L>\sin(\pi/(2n))$ then for all $t\in \bbR$,
$\#c(f^{-1}(t))<n$, \cite{M}. A similar scale also exists for co-Lipschitz
uniformly continuous maps \cite{M}.]
\end{rem}

\begin{rem}\lb{rem-n}
The estimate of the cardinality of the exceptional set $T_f$ is best possible,
in the sense that for any $n\in \bbN$ it is easy to construct examples of Lipschitz
quotient mappings  $f:
\bbR^2 \lra \bbR$ so that $\card(T_f)= n-1$ and $\#c(f^{-1}(t))=n$ for all $t\in \bbR\setminus
T_f$. In Figure~\ref{pic3} below, we present sketches of examples of such functions
for $n=2,3,4$. In each sketch, level sets for different values of $t$ are represented by
different styles of lines (within limits set by the drawing program ($\Xy$-pic)), and
the mapping $f$ is the distance in the $\ell_1$-metric from
the solid lines, which represent the preimage of 0,
 multiplied, in each component of the complement of the solid lines, by
the sign indicated.

\hspace{7mm}
  \xy <0.8cm,0cm>:
  (-4.5,2)*+{n=2},
  (-4.2,1.3)*+{T_f=\{0\}},
  (-3,1.7)*+{-},
(1.1,2.2)*+{0},
 (1.1,1.7)*+{-},
(-2.3,1.7)*+{+},
(-2.3,-1.7)*+{+},
(-3,2) ;
 (1,-2)**\dir{-},
 (-3,-2) ;
 (1,2)**\dir{-},
(3.6,2)*+{n=3},
  (4.2,1.3)*+{T_f=\{0,1\}},
(8,2) ;
 (8,-2)**\dir{-},
(6,2) ;
 (8,0)**\dir{-},
(6,-2) ;
 (8,0)**\dir{-},
(11,1) ;
 (10,0)**\dir{-},
(11,-1) ;
 (10,0)**\dir{-},
(9,2) ;
 (9,-2)**\dir{--},
(11,2) ;
 (9,0)**\dir{--},
(11,-2) ;
 (9,0)**\dir{--},
(6,1) ;
 (7,0)**\dir{--},
(6,-1) ;
 (7,0)**\dir{--},
(8,2.2)*+{0},
(9,2.2)*+{1},
(6,1.7)*+{+},
(6.7,1.7)*+{-},
(6.7,-1.7)*+{-},
(8.3,1.7)*+{+},
(11,.7)*+{-},
\endxy

\vspace{8mm}
 \hspace{3mm} \xy <.8cm,0cm>:
 (-2,1) ;
 (0,0)**\dir{-},
(-2,-1) ;
 (0,0)**\dir{-},
(2.4,3) ;
 (3.4,0)**\dir{-},
(2.4,-3) ;
 (3.4,0)**\dir{-},
(3.4,3) ;
 (3.4,-3)**\dir{-},
(-2,1.5) ;
 (1,0)**\dir{--},
(-2,-1.5) ;
 (1,0)**\dir{--},
(1.7,3) ;
 (2.7,0)**\dir{--},
(1.7,-3) ;
 (2.7,0)**\dir{--},
(-2,2) ;
 (2.05,0)**\dir{~},
(-2,-2) ;
 (2.05,0)**\dir{~},
(1,3) ;
 (2,-.05)**\dir{~},
(1,-3) ;
 (2,0.05)**\dir{~},
(4.1,3) ;
 (4.1,-3)**\dir{--},
(4.1,0) ;
 (7.1,-3)**\dir{--},
(4.1,0) ;
 (7.1,3)**\dir{--},
(4.8,0) ;
 (7.5,2.7)**\dir{-},
(4.8,0) ;
 (7.5,-2.7)**\dir{-},
(4.7,1.42) ;
 (4.7,3)**\dir{~},
(4.63,1.42) ;
 (6.2,3)**\dir{~},
(4.7,-1.5) ;
 (4.7,-3)**\dir{~},
(4.7,-1.7) ;
 (6.2,-3)**\dir{~},
(1.7,3.3)*+{1},
(.9,3.2)*+{2},
(2.4,3.3)*+{0},
(-1.7,1.1)*+{+},
(-1.7,.6)*+{-},
(2.4,2.3)*+{+},
(3,2.3)*+{-},
(3,-2.3)*+{-},
(3.7,2.3)*+{+},
(6.8,2.3)*+{+},
(7,1.8)*+{-},
(-5.3,3)*+{n=4},
  (-4.5,2.3)*+{T_f=\{0,1,2\}},
\endxy
\begin{fig}  \label{pic3}
\end{fig}
\end{rem}

For the proof of Theorem~\ref{main} we will need a large number of
auxiliary results concerning the number of components of level
sets $f^{-1}(t)$ and the end structure of boundaries of components
of the complements of $f^{-1}(t)$.  We start from a presentation
of these results and postpone the proof of Theorem~\ref{main} to
the end of this section.

Our first observation  relates the number of components of
$f^{-1}(t)$ with the number of components of $\bbR^{2} \setminus f^{-1}(t)$.

\begin{prop}\lb{number1}
Let $f:\bbR^2 \lra \bbR$ be a co-Lipschitz uniformly continuous map, $t \in \bbR$ and $K$ be
a connected
 component of the level set $f^{-1}(t)$.
Then
$$
\#c (\bbR^2 \setminus K)= \#e(K).
$$
In particular, $\#e(K)\ge 2$.

Moreover, each component of $\bbR \setminus K$ is homeomorphic with $\bbR^2$ and
its boundary is connected and has exactly 2 ends.
\end{prop}

\begin{proof}
By Theorem~\ref{graph}, $K$ is an unbounded finite graph, i.e.
 $K$ has a representation of the form:
$$
K=K_0 \cup \bigcup_{j=1}^n K_j,
$$
where $n = \#e(K), \ K_0$ is a compact connected tree, each $K_j$ is
a closed unbounded set homeomorphic with $[0,\infty)$, sets $\{K_j\}^n_{j=1}$ are
 mutually disjoint and $\card(K_j \cap K_0)=1$ for all $j$.

Let $\alpha (K_0)$ denote the number of vertices of $K_0$ and $\beta(K_0)$ the number of
 edges of $K_0$.  It is a basic fact of the graph theory
(see e.g. \cite[Theorem~IV.9 (attributed to Listing 1862)]{Konig}) that  since $K_0$
is a tree we have
$$
\alpha (K_0)= \beta(K_0)+1.
$$

Now let $\overline K$ be the closure of $K$ in the sphere $\mathcal{S}_2$ which is a
one-point compactification of $\bbR^2$.
Then $\overline K$ is a compact graph in $\mathcal{S}_2$, and, by the well-known
 Euler's formula, we have
$$
\alpha (\overline K)- \beta(\overline K) + \gamma (\overline K)=2,
$$
where $\alpha(\overline K)$ is the number of vertices of $\overline K, \beta (\overline
K)$  is the number of edges of $\overline K$, and $\gamma (\overline K)$ is the number of
components of $\mathcal{S}_2 \setminus \overline K$.  We have:
\begin{equation*}
\begin{split}
\alpha(\overline K)&= \alpha (K_0)+1,\\
\beta (\overline K)    &= \beta(K_0)+ n.
\end{split}
\end{equation*}

Thus
\begin{equation*}
\begin{split}
\#c(\bbR^2\setminus K)&= \gamma (\overline K) = 2 - \alpha (\overline K) + \beta (\overline K)\\
                     &= 2 - \alpha (K_0)-1 + \beta (K_0) + n\\
                     &= n = \#e (K),
\end{split}
\end{equation*}
as claimed (since the number of components of $\bbR^2 \setminus K$ is equal to  the
number of components of $S_2 \setminus \overline K$).

Since, by Corollary~\ref{sepplane}, $K$ separates the plane, thus
$\#e(K)=\#c(\bbR^2\setminus K)\ge2$.

Further, by \cite[Theorem~61.II.4]{Kur2}, the boundary of every component of  $\mathcal{S}_2
\setminus \overline K$ is a simple closed curve.
Thus by \cite[Theorem~61.V.1]{Kur2} each component of $\mathcal{S}_2 \setminus \overline K$,
and therefore also of $\bbR^2 \setminus K$, is homeomorphic with $\bbR^2$.  Since
$\infty$ belongs to the boundary of every component of $\mathcal{S}_2 \setminus \overline K$,
and since this boundary is a simple closed curve, we conclude that the order of
$\infty$, as a point of the boundary of any component of $\mathcal{S}_2 \setminus \overline
K$, is equal to two and thus this boundary has exactly 2 ends, and it is connected as a
subset of $\bbR^2$.
\end{proof}

Propostion~\ref{number1} has two useful consequences.

\begin{cor}\lb{number3}
Let $f: \bbR^2 \lra \bbR$ be a co-Lipschitz uniformly continuous map,  $t \in \bbR$
and  $A$ be a component of  $\bbR^2 \setminus f^{-1} (t)$. Then
$$
\#e (\Bd(A)) = 2 \# c (f^{-1}(t) \cap \Bd (A)).
$$
\end{cor}

\begin{proof}
Let $\{K_i\}^n_{i=1}$ be components of $f^{-1}(t)$.  If $\# c (f^{-1}(t)\cap \Bd
(A))=1$,  then $A$ is a component of, say, $\bbR^2 \setminus K_1$ and, by
Proposition~\ref{number1},
$$
\# e (\Bd(A))=2.
$$

We suppose, for the induction, that if $\#c (f^{-1}(t) \cap \Bd (A)) \leq k$, i.e.  if
$A$ is a component $\bbR^2 \setminus \bigcup^k_{j=1} K_j$, then
\begin{equation}\lb{count}
\# e (\Bd(A))= 2 \#c (f^{-1}(t) \cap \Bd (A)).
\end{equation}

Let $B$ be a component if $\bbR^2 \setminus \bigcup^{k+1}_{j=1} K_j$, say $B=A
\cap C$,  where $A$ is a component of $\bbR^2 \setminus \bigcup^k_{j=1} K_j$
and $C$ is a component of $\bbR^2 \setminus K_{k+1}$.

If $K_{k+1} \not \subset A$, then, by the connectedness of $K_{k+1}, K_{k+1} \cap A =
\emptyset$  and either $A \subset C$ or $A \cap C = \emptyset$.  Since $B \neq
\emptyset$, we obtain that $B=A$ and $B$ is a component of $\bbR^2 \setminus
\bigcup^k_{j=1}K_j$ and, by the inductive hypothesis, there is nothing to
prove.

Thus, without loss of generality, we assume that $K_{k+1} \subset A$ and
$\Bd(C) \subset K_{k+1} \subset A$.
Then
\begin{equation}\lb{bdry1}
\Bd(C)\subset \Bd (A \cap C).
\end{equation}

Similarly, if $\bigcup^k_{j=1}K_j \not \subset C$ then by the connectedness of
sets $\{K_j\}^{k}_{j=1}$, there exists $i_0 \leq k$ so that $K_{i_{0}} \not
\subset C$ and thus $K_{i_{0}} \cap C = \emptyset$.  Hence any component of $\bbR^2
\setminus K_{i_{0}}$ is either disjoint with $C$, or contains $C$.  Thus $B$ can be
represented as an intersection of components of
$\{\bbR^2 \setminus K_j\}^{k+1}_{j=1,
j \neq i_{0}}$ and, by the inductive hypothesis, we are done.

 Thus, without loss of
generality, we assume that $\bigcup^{k}_{j=1} K_j \subset C$.  Hence, as
before, $\Bd(A) \subset \bigcup^k_{j=1} K_j \subset C$ and
\begin{equation}\lb{bdry2}
\Bd(A)\subset \Bd (A \cap C).
\end{equation}

By \cite[Formula~6.II(8)]{Kur1}, we have
\begin{equation}\lb{bdry3}
\Bd(A \cap C) \subset \Bd (A) \cup \Bd (C).
\end{equation}

Combining \eqref{bdry1}, \eqref{bdry2} and \eqref{bdry3} we get
$$
\Bd (A \cap C)= \Bd(A) \cup
\Bd(C),
$$
and, since $\Bd (A)\cap \Bd(C)= \emptyset$, we conclude that
$$
\# e (\Bd(A \cap C))= \# e (\Bd(A))+ \# e (\Bd(C)).
$$

Thus, by \eqref{count} and \eqref{bdry1},
\begin{equation*}
\begin{split}
\# e (\Bd(B)) &= 2 \# c (f^{-1}(t) \cap \Bd(A)) +2 \\
             &= 2 \# c (f^{-1}(t) \cap(\Bd(A)\cup K_{k+1}))\\
             &= 2 \# c (f^{-1}(t) \cap \Bd (B)),
\end{split}
\end{equation*}
which ends the proof.
\end{proof}

\begin{cor}\lb{number2}
Let $f: \bbR^2 \lra \bbR$ be a co-Lipschitz uniformly continuous map and $t \in \bbR$.  Then
$$
\# c (f^{-1}(t)) + \# c(\bbR^2 \setminus f^{-1}(t)) = \# e(f^{-1}(t))+1.
$$
\end{cor}

\begin{proof}
By Theorems~\ref{finitecompl} and \ref{finite} both $\bbR^2 \setminus f^{-1}(t)$  and
$f^{-1}(t)$ have finite number of components.  Denote $ l = \# c (f^{-1}(t))$ and
let $\{K_j\}^l_{j=1}$ be the components of $f^{-1}(t)$.

By Proposition~\ref{number1},
$
\#c(\bbR^2 \setminus K_1) = \#e(K_1),
$
and $K_2$ is contained in exactly one of the components of $\bbR^2 \setminus K_1$, say
$C$.  Since $C$ is homeomorphic with $\bbR^2$, again by Proposition~\ref{number1} we
conclude that
$$
\# c (C \setminus K_2) = \# e(K_2),
$$
and thus
\begin{equation*}
\begin{split}
\#c (\bbR^2 \setminus (K_1 \cup K_2)) &= \#c (\bbR^2 \setminus K_1)-1
+ \#c (\bbR^2 \setminus K_2)\\
&= \#e(K_1)+ \#e (K_2)-1.
\end{split}
\end{equation*}
Proceeding by induction we get
\begin{equation*}
\begin{split}
\# c (\bbR^2 \setminus f^{-1}(t)) &= \# c
\Big(\bbR^2 \setminus (\bigcup^l_{j=1}K_j)\Big)\\
&= \sum^l_{j=1} \#e (K_j)- (l-1)\\
&= \# e (f^{-1}(t))+ 1 - \# c (f^{-1}(t)),
\end{split}
\end{equation*}
which ends the proof of the corollary.
\end{proof}

As a consequence of Theorems~\ref{finitecompl}, \ref{finite} and
Corollary~\ref{number2}  we immediately obtain:

\begin{cor}\lb{numberends}
Let $f:\bbR^2 \lra \bbR$ be a co-Lipschitz uniformly continuous map.  Then for every
$t \in \bbR$,  the number of ends of $f^{-1}(t)$ is finite and bounded by a function
depending only on
 the co-Lipschitz
 constant of $f$ and its modulus of uniform continuity.
\end{cor}

Our next goal is to show that the number of ends of $f^{-1}(t)$ is independent of $t$.
To achieve this we first prove that different ends of level sets $f^{-1}(t)$ are
``infinitely far away'' from each other, as on Figure~\ref{pic2} in the Introduction.
  To state this precisely we will use the
notation $d(X,Y)$ to denote the distance between sets $X,Y$ i.e.
$$
d(X,Y)\DEF \inf \{d(x,y): x \in X, y \in Y\}.
$$

\begin{prop}\lb{faraway}
Let $f: \bbR^2 \lra \bbR$ be a co-Lipschitz uniformly continuous map and
$t \in \bbR$.
If $f^{-1}(t)$
has $l$ components
$
\{K^{i}(t)\}_{i=1}^l
$,
and each component $K^i(t)$ has the following representation of   the form described in
Theorem~\ref{graph},
$$
K^{i}(t)=K^{i}_{0} (t)\cup \bigcup^{n(i)}_{j=1} K^{i}_{j}(t),
$$
then, for all $i_1, i_2 \in \{1,\cdots, l\}, j_1 \in \{1, \cdots
n (i_1)\}, j_2 \in \{1,\cdots, n(i_2)\}$, if the ordered pairs $(i_1, j_1),(i_2, j_2)$
are not the same, then
$$
\lim_{R \to \infty}d(K^{i_1}_{j_1}(t) \setminus B(0,R), K^{i_2}_{j_2}(t)
\setminus B(0,R))=\infty.
$$
\end{prop}

\begin{proof}
This result follows almost immediately from Lemma~\ref{bottleneck}. Without loss of
generality we assume that the co-Lipschitz  constant  of $f$ is $1$ and let $\Om(\cdot)$ be
the modulus of uniform continuity of $f$.
For any $R\in \bbR_+$ denote
$$
d_R =d(K^{i_1}_{j_1}(t) \setminus B(0,R),K^{i_2}_{j_2}(t) \setminus B(0,R)).
$$

Clearly $d_{R_1} \geq d_{R_2}$ when $R_1 \geq R_2$, thus, if  $\lim_{R \to
\infty} d_R \neq \infty$ then there exists $\alpha \in \bbR_+$ so that for all $R \in
\bbR$,
\begin{equation}\lb{small}
d_R \leq \al.
\end{equation}

Fix $x_1 \in K^{i_1}_{j_1}(t)$ and $y_1 \in K^{i_2}_{j_2} (t)$ so that
$$
d(x_1, y_1) \leq \al.
$$

Set
$$
\tilde{R}\DEF \| x_1 \| + 2\Om(\frac\al2)+4\al + 1.
$$

Then, by \eqref{small}, there exist $x_2 \in K^{i_1}_{j_1}(t) \setminus B(0,\tilde{R})$ and
$y_2 \in K^{i_2}_{j_2}(t)\setminus B(0,\tilde{R})$ with
$$
d(x_2, y_2)\leq \alpha.
$$

Since, for $\nu=1,2$, the sets $K^{i_\nu}_{j_\nu}(t)$ are connected
subsets of $f^{-1}(t)$, there exist  arcs $\sigma_\nu \subset
K^{i_\nu}_{j_\nu} (t)\subset f^{-1}(t)$ with endpoints $x_\nu,
y_\nu$.  Since $(i_1, j_1) \neq (i_2, j_2)$, the arcs $\sigma_\nu,
\nu=1,2$, are disjoint subcontinua of $f^{-1}(t)$, and hence, by
Lemma~\ref{bottleneck},
$$
d(x_1, x_2) \leq 2\Om(\frac\al2)+4\al.
$$
But
$$
d(x_1, x_2) \geq \big|\|x_2\|-\|x_1\|\big|\ge \tilde{R}-\|x_1\|=2\Om(\frac\al2)+4\al +1,
$$
and the resulting contradiction ends the proof of Proposition~\ref{faraway}.
\end{proof}

As an immediate corollary we obtain the following two facts which we state here for an easy
reference.

\begin{cor}\lb{distance}
Let $f: \bbR^2 \lra \bbR$ be a co-Lipschitz uniformly continuous map, $t \in \bbR$ and
$f^{-1}(t)$
have $l$  components
$
\{K^i (t)\}^l_{i=1}.
$
Then for all $i_1, i_2 \leq l$, $i_1 \neq i_2$, we have
$$
d (K^{i_1}(t), K^{i_2}(t))> 0.
$$
\end{cor}

\begin{proof}
This follows immediately from Proposition~\ref{faraway}.  Continuing the same notation
as above, let $R_0 \in \bbR$ be such that for all
$j_1 \in \{1, \cdots, n (i_1)\}, j_2 \in \{1,\cdots, n(i_2)\}$,
$$
d(K^{i_1}_{j_1}(t)\setminus B(0,R_0), K^{i_2}_{j_2}(t) \setminus B(0,R_0)) \geq 1.
$$

Since $K^{i_\nu}(t)\cap \overline{B(0,R_0)}$, for $\nu=1,2$, are compact
and disjoint we conclude that
$$
d(K^{i_1}(t), K^{i_2}(t)) \geq \min(1,d)>0,
$$
as desired.
\end{proof}

\begin{cor}\lb{farawayc}
Let $f: \bbR^2 \lra \bbR$ be a co-Lipschitz uniformly continuous map,
$t \in \bbR$, $R\in \bbR_+$
and $\{C_{R,i}\}_i$ are a collection of unbounded components of $f^{-1}(t)\setminus B(0,R)$.
Then, for all $i\ne j$,
$$
\lim_{R \to \infty}d(C_{R,i}, C_{R,j})=\infty.
$$
\end{cor}

\begin{rem}
After reading a preliminary version of this paper,  Maleva has strengthened the conclusion
of Corollary~\ref{farawayc}.  She proved \cite{M}, in the notation as above,
that there exists a constant
$\de>0$ depending only on the modulus of continuity of $f$ and its co-Lipschitz constant, so that
for every $t\in \bbR$ there exists $R(t)>0$ so that for all $R>R(t)$ and all $i\ne j$,
$$
d(C_{R,i}, C_{R,j})\ge \de R.
$$
This has consequences not only for the topology, but also for the allowable geometric
structure of $f^{-1}(t)$, e.g. $f^{-1}(t)$ cannot contain a parabola, see \cite{M}.
\end{rem}


As a consequence of Proposition~\ref{faraway}, we   obtain three  somewhat
technical facts which will be important for our further arguments.

\begin{lemma}\lb{endfacts}
Let $f:\bbR^2 \lra \bbR$ be a co-Lipschitz uniformly continuous map and $t_1, t_2 \in \bbR$.
 Then
\begin{itemize}
\item[(a)]
$\#e (f^{-1}(t_1)) = \#e (f^{-1}(t_2))$.
\item[(b)]
If $t_1>t_2$, $A$ is a component of $f^{-1}(t_2, \infty)$ or $A=\bbR^2$ and
 $f^{-1}(t_1, \infty) \cap A = \bigcup ^n_{\lambda=1} A_\lambda$, where $A_\lambda$
are components of $f^{-1}(t_1, \infty) \cap A $, then
$$
\sum^n_{\lambda=1} \#e (\Bd(A_\lambda)) = \#e (f^{-1}(t_1) \cap A).
$$
\item[(c)]
If $t_{1} > t_{2}$ and $A$ is a component of $f^{-1}(t_{2},\infty)$, then
$$
\#e(f^{-1}(t_{1}) \cap A) = \#e(\Bd(A)).
$$
\end{itemize}
\end{lemma}

\begin{proof}
As before we assume without loss of generality that the co-Lipschitz
constant  of $f$ is $1$.  By Corollary~\ref{number2} (or
a combination of earlier results in this paper) we know that both $f^{-1}(t_1)$ and
$f^{-1}(t_2)$ have a finite number of ends.  By Theorems~\ref{finite} and \ref{graph},
for $\nu=1,2, \ f^{-1}(t_\nu)$ can be represented as
\begin{equation}\lb{rep}
f^{-1}(t_\nu)=\bigcup^{l_\nu}_{i=1}\Big(K^i_0(t_\nu)\cup  \bigcup^{n_\nu(i)}_{j=1}
K^i_j(t_\nu)\Big),
\end{equation}
where $l_\nu, n_\nu(i)\in \bbN$ and $K_j^i(t_\nu)$ are mutually   disjoint unbounded rays in
$\bbR^2$, and sets $K^i_0(t_\nu)$ are compact.

Further, for $\nu=1,2$, the number of distinct rays $K^i_j(t_\nu)$  is finite and equals
the number of ends of $f^{-1}(t_\nu)$, i.e.
\begin{equation}\lb{card}
\card  \mathcal{R}(t_\nu) = \#e(f^{-1}(t_\nu)),
\end{equation}
where $\mathcal{R}(t_\nu)\DEF \{K^i_j (t_\nu): i=1,\dots, l_1, j=1, \dots,n_\nu (i)\}$.

By Proposition~\ref{faraway}, if $Y_1, Y_2$ are distinct rays of  the same level set
$f^{-1}(t_\nu)$, where $\nu\in\{1,2\}$, then
$$
\lim_{R \to \infty} d (Y_1 \setminus B(0,R), Y_2 \setminus B(0,R))= \infty.
$$

Thus, there exists $R_0 \in \bbR$ so that, for $\nu=1,2$, and for all $i \leq l_\nu$,
\begin{equation}\lb{bound1}
K^i_0 (t_\nu)\subset B(0, R_0-1),
\end{equation}
and so that, for any distinct  rays $Y_1, Y_2$ of the
same level set $f^{-1}(t_\nu)$, where $\nu=1$ or $2$,
\begin{equation}\lb{bound2}
d(Y_1 \setminus B(0,R_0) , Y_2 \setminus B(0,R_0)) \geq 1 + 4 | t_1 - t_2|.
\end{equation}

On the other hand, since $f$ is co-Lipschitz with constant $1$,  for every $x \in
f^{-1}(t_1)$ and for every $r > 0$,
$$
B(f(x),r)=(t_1-r,t_1 + r) \subset f(B(x,r)).
$$

Since $t_2 \in (t_1 - r, t_1 + r)$ when $r = 2 |t_1 - t_2|$, we  conclude that
\begin{equation}\lb{y}
\text{for  every} \ x \in f^{-1}(t_1)\ \text{ there  exists }
 \ y \in B(x, 2 |t_1 - t_2|) \cap f^{-1}(t_2).
\end{equation}

Now let $X \in \mathcal{R} (t_1)$  be a ray from
the representation of $f^{-1}(t_1)$ described in \eqref{rep}, and let $x \in X$ be
such that $\| x \| \geq R_0 + 2 |t_1 - t_2|$.  Then by \eqref{y} and \eqref{bound1}
there exists at least one ray $Y \in \mathcal{R} (t_2)$ so
that $d(x, Y \setminus B (0,R_0))< 2 |t_1 - t_2|$.

Suppose that there exist two distinct rays  $Y_1, Y_2 \in \mathcal{R} (t_2)$
so that for $\alpha =1,2$,
$$
d(x, Y_\alpha \setminus B(0,R_0))< 2 |t_1 - t_2|.
$$

But then we would have
$$
d(Y_1 \setminus B(0,R_0), Y_2\setminus B(0, R_0))< 4 | t_1-t_2|,
$$
which contradicts \eqref{bound2}.

 Thus we have described a one-to-one map $\gamma$
 from the
set of rays of $f^{-1}(t_1)$ into the set of rays of $f^{-1}(t_2)$, i.e. from
$\mathcal{R} (t_1)$ into $\mathcal{R} (t_2)$, and $\gamma$ operates in such a way that
for  every $X\in \mathcal{R} (t_1)$ and for every $x\in X$ with $\| x \| \geq R_0 + 2
|t_1 - t_2|$ we have
\begin{equation}\lb{gamma}
d(x, \gamma(X) \setminus B(0,R_0))< 2 |t_1 - t_2|.
\end{equation}

Since $\gamma$ is one-to-one,
  by \eqref{card} and by symmetry, we have
$$
\#e (f^{-1}(t_1))= \#e (f^{-1}(t_2)),
$$
which ends the proof of part (a).

Moreover, we conclude that $\gamma$ is   a bijection from $\mathcal{R} (t_1)$
onto $\mathcal{R} (t_2)$.

To prove part (b) we keep the same notation as above and we note that if
$
f^{-1}(t_1, \infty)\cap A= \bigcup^n_{\lambda =1} A_\lambda,
$
where $A_\lambda$ are components of $f^{-1}(t_1, \infty)\cap A$, then
\begin{equation}\lb{all}
f^{-1}(t_1)\cap A = \bigcup^n_{\lambda =1} \Bd (A_\lambda),
\end{equation}
and therefore
$$
\#e (f^{-1}(t_1)\cap A)=\#e \Big(\bigcup^n_{\lambda=1} \Bd (A_\lambda)\Big).
$$

Denote by $\mathcal{R}_A (t_1)$ the set of rays of $f^{-1}(t_1)$ contained in $A$, i.e.
$$\mathcal{R}_A (t_1)\DEF\{K^i_j (t_1)\subset\mathcal{R} (t_1): K^i_j (t_1)\subset A\}.$$

Since each ray $K^i_j (t_1)$ has exactly one end, we get
\begin{equation}\lb{card1}
\#e(f^{-1}(t_1)\cap A)= \card(\mathcal{R}_A (t_1)).
\end{equation}

By Theorem~\ref{graph}, for $R_0 \in \bbR_+$ defined above, and for each  $K^i_j (t_1)
\in\mathcal{R}_A (t_1)$, the set $K^i_j(t_1) \setminus B(0,R_0)$ has a unique
unbounded component; we will denote these components by $\{X_\al: \alpha =1, \dots,
\#e(f^{-1}(t_1)\cap A)\}$.  Note that, by Theorem~\ref{graph}, each $X_\alpha$ is
homeomorphic with $[0, \infty)$.
We will show that
\begin{equation}\lb{oneray}
\text{for each $\alpha \leq \#e(f^{-1}(t_1)\cap A)$ there exists a
unique $\lambda (\alpha)\leq n$ with} \ X_\al \subset \Bd (A_{\lambda(\alpha)}).
\end{equation}

Once \eqref{oneray} is established, part (b) follows easily.  Indeed,  by
\eqref{oneray} we can define sets
$$
E_\lambda  \DEF \{\alpha \leq \#e (f^{-1}(t_1)\cap A): X_\alpha \subset
\Bd (A_\lambda)\},
$$
and sets $E_\lambda$ are disjoint. Note that $\card(E_\lambda)=\#e(\Bd(A_\lambda))$.
 Moreover, by \eqref{all},
$$
\bigcup^n_{\lambda=1}\Bd (A_\lambda)  \supset \{X_\alpha: \alpha=1, \dots, \#e
(f^{-1}(t_1)\cap A) \},
$$
so $\bigcup^n_{\lambda=1}E_\lambda   = \{1, \dots, \#e
(f^{-1}(t_1)\cap A) \}$ and thus
$$
\#e (f^{-1}(t_1)\cap A)=
\sum^n_{\lambda=1} \card (E_\lambda)= \sum^n_{\lambda=1} \#e (\Bd(A_\lambda)),
$$
as desired.

To prove \eqref{oneray}, note that by  \eqref{all} and since sets $\Bd(A_\alpha)$ are
closed, for each $\alpha \leq \#e (f^{-1}(t_1)\cap A):$
$$
\overline{X_\alpha} \subset \bigcup^n_{\lambda=1} \Bd (A_\lambda).
$$

Thus
$$
\overline{X_\alpha}= \bigcup^n_{\lambda=1}(\overline{X_\alpha}\cap \Bd (A_\lambda)).
$$

Since $\overline{X_\al}$ is connected and sets $\overline{X_\alpha} \cap \Bd (A_\lambda)$
are closed, we conclude that either there exists a unique $\lambda(\alpha)$ so that,
for all $\lambda \neq \lambda (\alpha)$,
$$
\overline{X_\alpha} \cap \Bd (A_\lambda) = \emptyset,
$$
and in this case part (b) holds, or otherwise there exist $ \lambda_1, \lambda_2 \leq n,
\lambda_1 \neq \lambda_2$ so that
\begin{equation}\lb{alt}
(\overline{X_\alpha}\cap \Bd(A_{\lambda_{1}})) \cap (\overline{X_\alpha}\cap \Bd
(A_{\lambda_{2}})) \neq \emptyset.
\end{equation}

But this alternative leads to a contradiction.
Indeed, suppose that
$$
x \in \overline{X_\alpha} \cap \Bd (A_{\lambda_{1}}) \cap \Bd (A_{\lambda_{2}}).
$$

Since $\overline{X_\alpha}$ is a ray, i.e. a homeomorph of $[0,\infty)$,
 contained in one of the rays $\{K^i_j
(t_1)\}_{i,j}$ of $f^{-1} (t_1)\cap A$, and since, by \eqref{bound1},  $\overline{X_\alpha}$
is disjoint with all sets $ \{K^i_0 (t_1)\}_i$ we conclude that the order of the point $x$ in
$f^{-1}(t_1)$ is equal to $2$.  Hence, by Definition~\ref{deford}, for every
neighborhood $V$ of $x$, there exists a neighborhood of $x$ with $U \subset V$ and so
that $\card  (\Bd(U)\cap f^{-1}(t_1))=2$.  By \eqref{bound2} and since $f^{-1}(t_1)$ is
locally connected, we can choose $U \subset A$ so that $x \in U$, $\Bd (U)$ is a
simple curve, $\Bd(U)\cap f^{-1}(t_1)= \{x_1,x_2\}$ and $x$ belongs to an arc contained
in $f^{-1}(t_1)$ with endpoints $x_1$ and $x_2$.  Then, by the Theorem About The
$\theta$-Curve \cite[Theorem~61.II.2]{Kur2}, $U \setminus f^{-1}(t_1)$ has exactly two
components, and consequently $x$ belongs to the boundary of exactly two components of
$\bbR^2 \setminus f^{-1}(t_1)$.  Since $f$ is co-Lipschitz, it is not possible that
both of these components are contained in $f^{-1}(t_1, \infty)$,
which contradicts \eqref{alt} and  ends the proof of \eqref{oneray} and of part (b).

For part (c), let $\{A_{\nu}\}^{m}_{\nu=1}$  denote the
collection of components of
$f^{-1}(t_{2},\infty)$.  As above, let $\mathcal{R}_{A_{\nu}}(t_{1})$ denote the set of rays of
$f^{-1}(t_{1})$ contained in $A_{\nu}$.  By \eqref{card1}, for any $\nu \in \{1,\dots,m\}$,
\begin{equation} \lb{card 2}
\#e(f^{-1}(t_{1}) \cap A_{\nu})=\card (\mathcal{R}_{A_{\nu}}(t_{1})).
\end{equation}

Let $\gamma:\mathcal{R}(t_{1})\lra \mathcal{R}(t_{2})$ be the map defined in
part (a).  We will show that for all $\nu \in \{1,\dots,m\}$,
\begin{equation} \lb{gamma1}
X \in \mathcal{R}_{A_{\nu}}(t_{1})\ \ \ \Longrightarrow \ \ \ \gamma(X) \subset \Bd(A_{\nu}).
\end{equation}

Indeed, let $X \in \mathcal{R}_{A_{\nu}}(t_{1})$ and $x \in X$ with $\|x\|\geq R_{0} + 2
|t_{1}-t_{2}|$.  By \eqref{gamma}, there exists $y \in \gamma (X) \setminus
B(0,R_{0})$ so that $d(x,y)<2|t_{1}-t_{2}|$.  Let $I$ denote the interval with
endpoints $x$ and $y$.  If $y \notin \Bd(A_{\nu})$, then $I \cap \Bd(A_{\nu}) \not=
\emptyset$, since $\Bd(A_{\nu})$ separates between $x$ and $y$.  Thus there exists $y_{1}
\in I \cap \Bd(A_{\nu})$ so that $y_{1} \in f^{-1}(t_{2}),\|y_{1}\|\geq R_{0}$ and
$d(y,y_{1}) <2|t_{1}-t_{2}|$.  Thus, by \eqref{bound2}, $y_{1} \in \gamma(X)$ and
$\gamma(X) \cap \Bd(A_{\nu}) \not= \emptyset$.  Hence, by \eqref{oneray}, $\gamma(X)
\subset \Bd(A_{\nu})$ and \eqref{gamma1} holds.

Since $\gamma$ is one-to-one, \eqref{gamma1} immediately implies that
\begin{equation} \lb{card3}
\card (\mathcal{R}_{A_{\nu}}(t_{1})) \leq \#e(\Bd(A_{\nu})).
\end{equation}

Since sets $\{A_{\nu}\}^{m}_{\nu=1}$ are disjoint, by \eqref{card} and by part (b) applied
to $t_{2}$ and the set $A=\bbR^{2}$, we get
\begin{equation*}\begin{aligned}
\#e(f^{-1}(t_{1})) & = \card (\mathcal{R}(t_{1}))=\sum^{m}_{\nu=1} \card
(\mathcal{R}_{A_{\nu}}(t_{1}))\\
& \leq \sum^{m}_{\nu=1} \#e(\Bd(A_{\nu}))=\#e(f^{-1}(t_{2})).
\end{aligned}\end{equation*}

Since, by part (a), $\#e(f^{-1}(t_{1})) = \#e(f^{-1}(t_{2}))$, we conclude that, for
all $\nu \in \{1,\dots,m\}$,
$$ \#e (f^{-1}(t_{1})\cap A_{\nu})=\card (\mathcal{R}_{A_{\nu}}(t_{1}))
= \#e (\Bd(A_{\nu})),$$
 which ends the proof of part (c).
\end{proof}

For the proof of the main theorem we will need one more lemma.

\begin{lemma} \lb{comp}
Let $f:\bbR^{2} \lra \bbR$ be a co-Lipschitz uniformly continuous map, $t \in \bbR$ and $A$ be
a component of $f^{-1}(t,\infty)$.  Then:
\begin{itemize}
\item[(a)]
for every $s > t$,
$$
\#c(f^{-1}(s)\cap A) \leq \#c(f^{-1}(t) \cap \Bd(A));
$$
\item[(b)]
there exists $\e  > 0$ so that for every $s \in (t,t+\e)$,
$$
\#c(f^{-1}(s) \cap A) = \#c(f^{-1}(t) \cap \Bd(A));
$$
\item[(c)] let $t_{1}>t$ and let $\{C_{\nu}\}^{k}_{\nu=1}$ be components of
$f^{-1}(t_{1},\infty)$ which are contained in $A$,  then
$$\sum^{k}_{\nu=1}
\#c(f^{-1}(t_{1}) \cap \Bd (C_{\nu}))=\#c(f^{-1}(t)\cap \Bd (A));$$
\item[(d)] let
$m=\#c(f^{-1}(t) \cap \Bd(A))$,  then there exists a set $T_{A} \subset (t,\infty)$
with $\card (T_{A}) \leq m-1$, so that for every $s \in (t,\infty) \setminus T_{A}$,
$$
\#c(f^{-1}(s) \cap A) =\#c(f^{-1}(t) \cap \Bd(A)).
$$
\end{itemize}
\end{lemma}

\begin{proof}
Let $s>t$.  By Lemma~\ref{endfacts}(b) and (c) and by Corollary~\ref{number3},
$$\#e(f^{-1}(s)\cap A) = \#e(\Bd(A))=2 \#c(f^{-1}(t) \cap \Bd(A)).$$

Since, by Proposition~\ref{number1}, each component of $f^{-1}(s) \cap A$ has at
least 2 ends, we get that
$$\#c(f^{-1}(s) \cap A) \leq \frac{1}{2} \#e(f^{-1}(s)
\cap A) = \#c(f^{-1}(t) \cap \Bd(A)),$$
 and part (a) is proven.

For part (b) we assume, without loss of generality, that the co-Lipschitz constant of $f$
is equal to $1$.  Denote
$$
l \DEF \#c(f^{-1}(t) \cap \Bd (A)),
$$
and let $\{L_{j}\}^{l}_{j=1}$ be   the components of $f^{-1}(t) \cap \Bd(A)$.

By Proposition~\ref{number1}, the intersection of any component of $f^{-1}(t)$ with
$\Bd(A)$ is connected and therefore each $L_{j}$ is contained in a different component
of $f^{-1}(t)$.  Hence, by Corollary~\ref{distance}, there exists $\delta > 0$ so that,
for all $i,j \in \{1,\dots,l\},i \not= j,$
$$
d(L_{i},L_{j}) \geq \delta.
$$

Define for $j \in \{1,\dots,l\}$,
$$
U_{j} \DEF \bigcup_{x \in L_{j}} B(x,\frac{\delta}{2}) \cap A.
$$
Then $\{U_{j}\}^{l}_{j=1}$ are connected, mutually disjoint, open subsets of $A$.

 We
note that, for any $y \in \bbR^{2}$,
\begin{equation} \lb{big}
d(y,f^{-1}(t)) \geq \frac{\delta}{2}\ \ \ \Longrightarrow \ \ \ |f(y)-t| \geq \frac{\delta}{2}.
\end{equation}

Indeed, if $d(y,f^{-1}(t))\geq \frac{\delta}{2}$, then $t \notin
f(B(y,\frac{\delta}{2}))$.  But, since $f$ is co-Lipschitz with constant 1,
$f(B(y,\frac{\delta}{2}))\supset B(f(y),\frac{\delta}{2})$. Thus $t \notin
B(f(y),\frac{\delta}{2})$ and \eqref{big} holds.
Hence, for all $s \in (t, t+\frac{\delta}{2})$,
\begin{equation} \lb{u}
f^{-1}(s) \cap A \subset \bigcup^{l}_{j=1} U_{j}.
\end{equation}

Now fix $j_{0} \in \{1,\dots,l\}$, and let $x \in L_{j_{0}}$ and   $y \in \Bd
(U_{j_{0}}) \setminus \Bd(A)$. Then $d(y,f^{-1}(t)) \geq \frac{\delta}{2}$.

Moreover, since $L_{j_{0}}$ is locally connected and $U_{j_{0}}$ is open,
$\overline{U_{j_{0}}}$ is arcwise connected and there exists a continuous function
$\sigma:[0,1]\lra\overline{U_{j_{0}}}$ so that $\sigma(0)=x, \sigma(1)=y$ and
$\sigma(\lambda) \in
U_{j_{0}}$ for $\lambda \in (0,1)$.  Define $g:[0,1]\lra\bbR$ as $g=f \circ \sigma$.
Then
\begin{equation*}\begin{aligned}
g(0) &=f(x)=t,\\
g(1) &=f(y) \geq t+\frac{\delta}{2},\ \  \text{ by \eqref{big}.}
\end{aligned}
\end{equation*}

By the Intermediate Value Theorem, for every $s \in (t,t+\frac{\delta}{2})$, there
exists at least one $\lambda_{s} \in (0,1)$ so that
$$s=g(\lambda_{s})=f(\sigma(\lambda_{s})).$$

Since $\sigma(\lambda_{s}) \in U_{j_{0}}$, we conclude that, for every $j_{0} \in
\{1,\dots,l\}$ and every $s \in (t,t+\frac{\delta}{2})$,
$$
f^{-1}(s) \cap U_{j_{0}} \not= \emptyset.
$$

Since sets $\{U_{j}\}^{l}_{j=1}$ are mutually disjoint and by \eqref{u}, we get
that for all $s \in (t,t+\frac{\delta}{2})$,
$$
\#c(f^{-1}(s) \cap A) \geq l = \#c(f^{-1}(t) \cap \Bd(A)).
$$
This, together with part (a), concludes the proof of part (b).

Part (c) follows by the following computation:
\begin{equation*}
\begin{alignedat}{5}
\sum^{k}_{\nu=1} \#c(f^{-1}(t_{1}) \cap \Bd (C_{\nu}))&=\frac{1}{2} \sum^{k}_{\nu=1}
\#e (\Bd (C_{\nu})), &\qquad   & \quad\text{by Corollary~\ref{number3}},\\
&=\frac{1}{2} \#e(f^{-1}(t_{1}) \cap A),&   & \quad\text{by Lemma~\ref{endfacts}(b)},\\
&=\frac{1}{2} \#e(\Bd (A)),& & \quad\text{by
 Lemma~\ref{endfacts}(c)},\\
&=\#c(f^{-1}(t) \cap \Bd (A)),& & \quad \text{by  Corollary~\ref{number3}.}
\end{alignedat}
\end{equation*}

To prove part (d), we proceed inductively with respect to $m$.

 If
$m=1$, then by part (a), for every $s \in (t,\infty)$,
$$
\#c(f^{-1}(s) \cap A) \leq 1,
$$
and since $f(A)=(t,\infty)$, we have $f^{-1}(s) \cap A \not=
\emptyset$, and hence $\#c(f^{-1}(s) \cap A) \geq 1$.  Therefore
part (d) holds with $T_{A}=\emptyset$, as desired.

 For the
induction, we assume that part (d) holds for all   $m<m_{0}$, where
$m_{0}\geq 2$.

Now suppose that
$$
\#c(f^{-1}(t) \cap \Bd (A))=m_{0} \geq 2.
$$
Define
$$t_{1} \DEF \sup \{\tau \in (t,\infty):\forall s \in
(t,\tau)\ \  \#c(f^{-1}(s) \cap A)=\#c(f^{-1}(t) \cap \Bd(A))\}.$$

If $t_{1}=\infty$ there is nothing to prove, so suppose that
$t_{1} < \infty$.  By part (b), $t_{1}>t$ and
$$\#c(f^{-1}(t_{1})
\cap A) \not= \#c(f^{-1}(t) \cap \Bd(A)).$$

By part (a), this implies that
\begin{equation} \lb{few}
\#c(f^{-1}(t_{1}) \cap A) < \#c(f^{-1}(t) \cap \Bd(A)).
\end{equation}

Let $\{C_{\nu}\}^{k}_{\nu=1}$ denote all components of
$f^{-1}(t_{1},\infty)\cap A$.
Then, by \eqref{few}, for each $\nu \leq k$,
\begin{equation} \begin{aligned} \lb{few2}
\#c(f^{-1}(t_{1}) \cap \Bd (C_{\nu})) &\leq \#c(f^{-1}(t_{1}) \cap
A)\\
&< \#c(f^{-1}(t) \cap \Bd(A))=m_{0}.
\end{aligned} \end{equation}

Hence, by the inductive hypothesis, for each $\nu \leq k$ there exists a set
$T_{\nu}=T_{C_{\nu}} \subset (t_{1},\infty)$ with $\card(T_{\nu}) \leq
\#c(f^{-1}(t_{1}) \cap \Bd(C_{\nu}))-1$, so that for every $s \in (t_{1},\infty)
\setminus T_{\nu}$,
\begin{equation} \lb{s}
\#(f^{-1}(s) \cap C_{\nu})=\#c(f^{-1}(t_{1})\cap \Bd(C_{\nu})).
\end{equation}

Set
$$
T_{A}=\bigcup^{k}_{\nu=1} T_{\nu} \cup \{t_{1}\}.
$$
Then for every $s \in (t_{1},\infty) \setminus T_{A}$, we have:
\begin{equation*}
\begin{alignedat}{5}
\#c(f^{-1}(s) \cap A)&=\sum^{k}_{\nu=1} \#c (f^{-1}(s) \cap C_{\nu})&\qquad &  \quad\\
&= \sum^{k}_{\nu=1} \#c(f^{-1}(t_{1}) \cap \Bd(C_{\nu})),& & \text{by
 \eqref{s}},\\
&= \#c(f^{-1}(t) \cap \Bd(A)),& & \text{by   part (c).}
\end{alignedat}
 \end{equation*}

To finish the proof we only need to estimate the cardinality of the set $T_{A}$.
We have
\begin{equation*}
\begin{alignedat}{5}
\card (T_{A}) &\leq \sum^{k}_{\nu=1} \card (T_{\nu}) +1  &\quad &\\
&\leq \sum^{k}_{\nu=1} \big[\#c(f^{-1}(t_{1}) \cap \Bd(C_{\nu}))-1\big]+1& &\\
&= \#c(f^{-1}(t) \cap \Bd(A)) + (1-k), & &\text{by  part (c).}\\
&\leq \#c(f^{-1}(t) \cap \Bd(A))-1,& &\text{since, by part (c) and \eqref{few2}, $k \geq
2$},
\end{alignedat}
\end{equation*}
which ends the proof of part (d).
\end{proof}

We are now ready for the proof of our main theorem.

\begin{proof}[Proof of Theorem~\ref{main}]
To prove part (1) we note that by Corollary~\ref{numberends},
number of ends of any level set $f^{-1}(t)$, for $t \in \bbR$, is finite
and bounded by a constant $M$ depending only on the
co-Lipschitz constant of $f$ and its modulus of uniform continuity.

By Lemma~\ref{endfacts}(a), $\#e(f^{-1}(t))$ does not depend on
the value of $t \in \bbR$.  To see that $\#e(f^{-1}(t))$ is even,
let $\{A_{\lambda}(t)\}^{l}_{\lambda=1}$ be the components of
$\bbR^{2} \setminus f^{-1}(t)$.  Then we have
\begin{equation*}
\begin{alignedat}{5}
\#e(f^{-1}(t))&= \sum^{l}_{\lambda=1}\#e(\Bd(A_{\lambda}(t))),&\quad& \text{by
 Lemma~\ref{endfacts}(b)},\\
&=2\sum^{l}_{\lambda=1} \#c(f^{-1}(t) \cap \Bd(A_{\lambda}(t))),& &
\text{by  Corollary~\ref{number3}}.
\end{alignedat}
\end{equation*}

Thus $\#e(f^{-1}(t))$ is even.

The moreover statement follows from Corollary~\ref{farawayc},
 and hence part (1) is proven.

For the proof of part (2), let $t_{0}$ be any real number, say
$t_{0}=0$, and let $\{A_{\nu}\}^{l}_{\nu=1}$ be all the components of
$f^{-1}(0,\infty)$.
By Lemma~\ref{comp}(d), for every $\nu \leq l$, there exists a set $T_{A_{\nu}}
\subset (0,\infty)$ with $\card (T_{A_{\nu}}) \leq \#c(f^{-1}(0) \cap
\Bd(A_{\nu}))-1$, so that for all $s \in (0,\infty) \setminus T_{A_{\nu}}$,
\begin{equation} \lb{last2}
\#c(f^{-1}(s) \cap A_{\nu})=\#c(f^{-1}(0) \cap \Bd(A_{\nu})).
\end{equation}

Define
$$
T^0_{+} = \bigcup^{l}_{\nu=1} T_{A_{\nu}}.
$$
Note that
\begin{equation}\lb{last1}
\begin{alignedat}{5}
\sum^{l}_{\nu=1} \#c(f^{-1}(0) \cap \Bd (A_{\nu})) &= \frac{1}{2}
\sum^{l}_{\nu=1} \#e(\Bd(A_{\nu})),&\qquad& \text{by   Corollary~\ref{number3}},\\
&= \frac{1}{2} \#e(f^{-1}(0)),& & \text{by   Lemma~\ref{endfacts}(b)}.
\end{alignedat}
\end{equation}

Therefore for every $s \in (0,\infty) \setminus T^0_{+}$, we have:
\begin{equation*}
\begin{alignedat}{5}
\#c(f^{-1}(s)) &= \sum^{l}_{\nu=1} \#c(f^{-1}(s) \cap A_{\nu})& &\\
&=\sum^{l}_{\nu=1} \#c(f^{-1}(0) \cap \Bd(A_{\nu})),&\quad& \text{by   \eqref{last2}},\\
&= \frac{1}{2} \#e(f^{-1}(0)),& & \text{by \eqref{last1}},\\
&= n,& & \text{by   part (1)}.
\end{alignedat}
\end{equation*}
Similarly
\begin{equation*}
 \begin{alignedat}{5}
\card (T^0_{+}) &\leq \sum^{l}_{\nu=1} \card (T_{A_{\nu}})
\leq \sum^{l}_{\nu=1}\big[\#c(f^{-1}(0) \cap \Bd (A_{\nu}))-1\big]& &\\
&= \frac{1}{2} \#e(f^{-1}(0))-l,& & \text{by   \eqref{last1}},\\
&\leq n-1,& & \text{by   part (1)}.
\end{alignedat}
\end{equation*}

Next, we note that in an identical way (e.g. by replacing function $f$ by $-f$) one can
define a set $T^0_-\subset(-\infty,0)$ with $\card  (T^0_-)\leq n-1$, so that for every $s \in
(-\infty,0) \setminus T^0_-$,
$
\#c(f^{-1}(s)) = n.
$
Define
$$T^0_{f}=T^0_{+} \cup T^0_- \cup \{0\}.$$

Clearly, $\card (T^0_{f}) \leq 2n-1$.
Now, let $t_{00}\in\bbR$ be such that $t_{00}<t$ for every $t\in T_f^0$.
Then we clearly have $\#c(f^{-1}(s)) = n$, for every $s\le t_{00}$. Same way as was done
above for $t_0=0$, we construct a set $T_+^{00}\subset(t_{00},\infty)$ so that for every
$s \in(t_{00},\infty) \setminus T^{00}_+$, we have
$\#c(f^{-1}(s)) = n$ and $ \card(T_+^{00})\le n-1$. Then part (2a) holds for
$T_f\DEF T_+^{00}$.

Part (2b) follows immediately by Corollary~\ref{number2}.

For part (2c), note that by Proposition~\ref{number1}, each
component of $f^{-1}(t)$ has at least 2 ends.  Since by parts (1)
and (2a) for all $t \in \bbR \setminus T_{f},\ f^{-1}(t)$  has $n$
components and $2n$ ends, we conclude that each component of
$f^{-1}(t)$ has exactly 2 ends.  Hence, by
Proposition~\ref{number1}, each component $K$ of $f^{-1}(t)$
separates the plane into exactly 2 components.  Further, by
Theorem~\ref{graph}, $K$ has a representation of the form
$$
K=K_{0} \cup K_{1} \cup K_{2},
$$
where $K_{1} \cup K_{2} = \emptyset$, and for $i=1,2$, $K_{i} \cap K_{0}$
consists of exactly one point which is on endpoint of both $K_{i}$
and $K_{0}$, $K_{1}$ and $K_{2}$ are both homeomorphic with
$[0,\infty)$, and $K_{0}$ is a compact connected tree with exactly
2 endpoints.  Thus $K_{0}$ is homeomorphic with $[0,1]$ (either by
the construction of $K_{0}$ described in the proof of
Theorem~\ref{graph}, or by the classical characterization of the
interval as a continuum with exactly 2 non-cut points cf. e.g.
\cite[Theorem~III.(6.2)]{WhyAT}). Therefore $K$ is homeomorphic with
$(-\infty,0] \cup [0,1] \cup [1,\infty) = (-\infty,\infty)$, which
ends the proof of part (2c).

For part (3a) we note that for all $t_{i} \in T_{f}$,
$\#c(f^{-1}(t_{i})) \not= n$.  Since $\#e(f^{-1}(t_{i}))=2n$ by
part (1), and each component has at least 2 ends, by
Proposition~\ref{number1}, we conclude that
$\#c(f^{-1}(t_{i}))<n$, i.e. part (3a) holds.  Parts (3b) and (3c)
follow immediately from Corollary~\ref{number2} and
Theorem~\ref{graph}, respectively.
\end{proof}

\section{Example
}
\lb{example}

In this section  we present an example, mentioned in the Introduction,
 of a uniform quotient map  $f:
\bbR^2 \lra \bbR$ with a non-locally connected level set.

\begin{example} \lb{exnlc}
There exists a uniform quotient map $f:\bbR^2 \lra \bbR$ so that $f^{-1}(0)$  is not
locally connected.
\end{example}

\begin{proof}[Construction]
Let $z_n = (\frac{1}{n},(-1)^n)\in \bbR^2$ for $n \in \bbZ \setminus
 \{0\}$, and let $I_n$ be a segment in $\bbR^2$ with endpoints $z_n,\ z_{n+1}$,   when $n>0$,
 or $z_n,\ z_{n-1}$ when $n<0$.  Let $I_0$ be the vertical segment with endpoints
$(0,1)$ and $(0,-1)$, and let $I_+, I_-$ be the following two half-lines:
$$
I_+=\{(x,-1): x \geq1\},\ \ \ I_- = \{(x,-1): x \leq-1\}
$$

Define $K$ to be the sum of all these segments
$$
K \overset {\text{def}}{=} \bigcup_{n \in \bbZ}
I_n \cup I_+ \cup I_-
$$

\vspace{-1mm}

  \xy <2cm,0cm>:
 (-6.2,0) ;
 (-5,0)**@{-},
 (-5,0) ;
 (-4.5,1)**@{-},
 (-4.5,1) ;
 (-4,0)**@{-},
 (-4,0) ;
 (-3.7,1)**@{-},
 (-3.4,0) ;
 (-3.7,1)**@{-},
 (-3.4,0) ;
 (-3.2,1)**@{-},
 (-2.8,1) *={\dots},
 (-2.8,0) *={\dots},
 (-2.4,0);
 (-2.4,1) **@{-},
(1.4,0) ;
 (0.2,0)**@{-},
 (0.1,1.2)*=0{+},
(0.1,-0.2)*=0{-},
 (0.2,0) ;
 (-0.3,1)**@{-},
 (-0.3,1) ;
 (-0.8,0)**@{-},
 (-0.8,0) ;
 (-1.1,1)**@{-},
 (-1.4,0) ;
 (-1.1,1)**@{-},
 (-1.4,0) ;
 (-1.6,1)**@{-},
 (-2.0,1) *={\dots},
 (-2.0,0) *={\dots},
\endxy
\vspace{-1mm}
\begin{fig} { Set $K$.}\label{f1}
\end{fig}

Set $K$ is connected but not locally connected and it separates the plane into two
regions.   We define the map $f:\bbR^2 \lra \bbR$ as the distance from $K$ multiplied
in each component of $K$ by the sign indicated.
Then, clearly, $K=f^{-1}(0)$.  It is also clear that $f$ is Lipschitz.
We will show that $f$ is co-uniformly continuous with
\begin{equation}\lb{omega}
\omega(r)=
\begin{cases}
\frac{r^3}{16000}          &\text{if $r < \frac{1}{10}$,} \\
\frac{1}{16 \cdot 10^6}    &\text{if $r \geq \frac{1}{10}.$}
\end{cases}
\end{equation}

We achieve this in a number of steps.
\begin{step}\lb{1}
Let $A_n = (\frac{1}{n},0),\ B_n = (\frac{1}{n},(-1)^n),\ C_n
=(\frac{1}{n-1},(-1)^{n-1}),\ D_n = (\frac{1}{n+1},(-1)^{n+1})$ for $n \in \bbZ
\setminus \{0,1,-1\}$. Let $\al_n =\sphericalangle A_n B_n C_n$ and
$\be_n = \sphericalangle A_n B_n D_n$, where
both  angles are assumed to be positive.  Then
\begin{equation*}
\begin{split}
\frac{1}{2|n|(|n|-1)} &\geq \sin \al_n \geq  \frac{1}{3|n|(|n|-1)}\\
\frac{1}{2|n|(|n|+1)} &\geq \sin \be_n \geq  \frac{1}{3|n|(|n|+1)}.
\end{split}
\end{equation*}
\end{step}

\begin{proof}[Proof of Step~\ref{1}]
Without loss of generality, we assume that $n >1$.   We illustrate $\al_n$ and $\be_n$
on Figure~\ref{f2}.

\hspace{35mm}  \xy <0.8cm,0cm>:
 (-3.8,3) ;
 (0,-3)**@{-},
 (0,-3) ;
 (4,3)**@{-},
 (-5,0) ;
 (5,0)**@{-},
 (-3.8,0) ;
 (-3.8,3)**@{.},
 (0,0) ;
 (0,-3)**@{.},
 (4,0) ;
 (4,3)**@{.},
(4.1,3.2)*+{C_n},
(-3.7,3.2)*+{D_n},
(0,-3.3)*+{B_n},
(-3.8,-0.4)*+{\frac{1}{n+1}},
(4,-0.4)*+{\frac{1}{n-1}},
(2,0)*+\dir{*},
(0,0)*+\dir{*},
 (1.5,0.1); (0.1,0.1) **\frm{^\}},
 (1.5,0.5) *+{\frac{1}{2n(n-1)}},
 (2.3,-0.5) *+{\frac{2n-1}{2n(n-1)}},
 (-0.1,0.5)*+{\frac1n},
 (-0.2,-0.4)*+{A_n},
 (0.3,-2.1)*+{\alpha_n},
 (-0.2,-2.1)*+{\beta_n},
 \endxy
\begin{fig} \label{f2}
\end{fig}

It is not difficult to compute that, as indicated on Figure~\ref{f2},
$$
\sin \al_n = \frac{\frac{1}{2n(n-1)}}{\sqrt{1+(\frac{1}{2n(n-1)})^2}}
= \frac{1}{\sqrt{4n^2 (n-1)^2+1}}.
$$
Thus,
$$
\frac{1}{2n(n-1)} = \frac{1}{\sqrt{4n^2 (n-1)^2}}
\geq \sin \al_n \geq \frac{1}{\sqrt{9n^2 (n-1)^2}}
= \frac{1}{3n (n-1)}.
$$
Similarly,
$$
\sin \be_n = \frac{\frac{1}{2n(n+1)}}{\sqrt{1+(\frac{1}{2n(n+1)})^2}}
= \frac{1}{\sqrt{4n^2 (n+1)^2+1}}.
$$
Thus
$$
\frac{1}{2n(n+1)} = \frac{1}{\sqrt{4n^2 (n+1)^2}}
\geq \sin \be_n \geq \frac{1}{\sqrt{9n^2 (n+1)^2}}
= \frac{1}{3n (n+1)}.
$$
\end{proof}

\begin{step}\lb{3}
If $x=(x_1, x_2)\in K$ and $x_1=0$ then for $r \leq \frac{1}{10}$
\begin{equation}\lb{rcube}
f(B(x,r)) \supset B (f(x), \frac{r^3}{2000}).
\end{equation}
\end{step}

\begin{proof}[Proof of Step~\ref{3}]
We first consider the case when $|x_2| \leq 1$, as illustrated on Figure~\ref{f5}.

Let $n$ be the smallest odd number so that
\begin{equation}\lb{n}
\frac{1}{n-1} \leq \frac{r}{32}.
\end{equation}


  \hspace{60mm}\xy <1.4cm,0cm>:
 (-1.7,0);
 (-1.7,2) **@{-},
 (0.2,0) ;
 (-0.3,2)**@{-},
 (-0.3,2) ;
 (-0.8,0)**@{-},
 (-0.8,0) ;
 (-1.1,2)**@{-},
 (-1.4,0) ;
 (-1.1,2)**@{-},
 (-1.4,0) ;
 (-1.6,2)**@{-},
(-1.7,0.2)*+{\cdot},
(-1.9,0.2)*+{x},
 (-1.7,0.2)*\cir<8mm>{}="*",
 (-1.4,0);
 (-1.4,0.55);
(-1.7,0.2)**@{-},
(-1.4,0.55);
(-1.4,0)**@{.},
(-1.4,0.55)*{\cdot},
(-1.35,0.7)*{y},
(-1.7,0.2);
(-1.25,0.2)**@{.},
(-1.4,0.35); (-1.6, -.4) **\crv{ (-1.8,0.1)&(-1.6,0)} ?<(0)*\dir{<},
(-1.6,-.5)*{\beta_n},
\endxy

\begin{fig} { }\label{f5}
\end{fig}

Then there exists $y=(y_1,y_2) \in B (x,r)$ so that  $y_1= \frac{1}{n}$ and $y_2 \geq
x_2 + \frac{r}{2}.$  Then
$$
f(y)= d(y,K) \geq \frac{r}{2} \sin \be_n \geq \frac{r}{2} \cdot \frac{1}{3n(n+1)}.
$$
By \eqref{n} and since $r \leq \frac{1}{10}$ we see that $n > 4$  and $\frac{1}{n-3} >
\frac{r}{2}$. Thus
\begin{equation*}
\begin{alignedat}{3}
\frac{1}{n}   &=\frac{n-3}{n} \cdot \frac{1}{n-3} &\ \geq \frac{1}{4} \cdot \frac{r}{2},\\
\frac{1}{n+1} &= \frac{n-3}{n+1}\cdot \frac{1}{n-3} &\ \geq \frac{1}{5}\cdot \frac{r}{2}.
\end{alignedat}
\end{equation*}
Hence
$$
f(y)= d(y,k) \geq \frac{r}{2} \cdot \frac{1}{3} \cdot \frac{1}{4}
\frac {r}{2} \cdot \frac{1}{5} \cdot \frac{r}{2} = \frac{r^3}{480}.
$$

Similarly there exists $z=(z_1, z_2) \in B(x,r)$ so that $z_1 = \frac{1}{n+1}$  and
$z_2 \leq x_2 - \frac{r}{2}$.  Then $f(z)= -d(z, K)$ and
$$
d(Z, K)\geq \frac{r}{2} \sin \be_{n+1} \geq \frac{r}{2}\cdot  \frac{1}{3(n+1)(n+2)}  \geq
\frac{r}{2} \cdot \frac{1}{3} \cdot \frac{1}{5}\cdot  \frac{r}{2} \cdot
\frac{1}{6}\frac{r}{2} = \frac {r^3}{720}.
$$
Thus \eqref{rcube} is satisfied.
\end{proof}

\begin{step}\lb{4}
If $x=(x_1, x_2) \in K$ and $x_1 \neq 0$ then for $r \leq \frac{1}{10}$,
 \eqref{rcube} is satisfied.
\end{step}

\begin{proof}[Proof of Step~\ref{4}]
Since $x_1 \neq 0$, thus there exists $m \in \bbZ$ so that $x$ belongs to $I_m$, the segment
which connects points $(\frac{1}{m},(-1)^m)$ and $(\frac{1}{m+1},(-1)^{m+1})$.
Without loss of generality we will assume that $m>0$.  Let $n$
denote the smallest odd number so that
$$\frac{1}{n-1}\leq \frac{r}{2}.$$
If $n < m-1$
we proceed in a way very similar to Step~\ref{3}.  Since $x_1 > 0$, we see that there
exists $y=(y_1, y_2) \in B (x,r)$ and $z=(z_1, z_2)\in B(x,r)$ so that $y_1 =
\frac{1}{n},\ y_2 \geq x_2 + \frac{r}{2},\ z_1=\frac{1}{n+1},\ z_2 \leq x_2 -
\frac{r}{2}$.  Then $f(y)= d(y,K)$ and $f(z)= -d (z, K)$.
Further, similarly  as in Step~\ref{3},
\begin{equation*}
\begin{split}
d (y,K) &\geq \frac{r^3}{480},\\
d (z,K) &\geq \frac{r^3}{720}.
\end{split}
\end{equation*}
Thus \eqref{rcube} is satisfied.

If $n \geq m-1$ and $m>3$ (the case when $m \leq 3$ is done similarly and  we
leave the details to the interested reader) then
\begin{equation} \lb{m}
\frac{1}{m-3}> \frac{r}{2}.
\end{equation}

Now let $t \in [0,1]$ be such that  $x=t(\frac{1}{m}, (-1)^m)+
(1-t)(\frac{1}{m+1},(-1)^{m+1})$, see Figure~\ref{f7}.


\hspace{27mm}  \xy <0.7cm,0cm>:
 (-3.6,3) ;
 (0,-3)**@{-},
 (0,-3) ;
 (3.8,3)**@{-},
 (0,3) ;
 (0,-3)**@{.},
 (-1.2,-1)*+\dir{*},
 (-1.5,-1)*+{x},
 (-.75,-.35)*+{d},
 (0,-.28)*+\dir{*},
 (0,-.28);
 (-1.2,-1)**@{-},
 (-1,-.88);
 (-.88,-1.08)**@{-},
 (-1.08,-1.2);
 (-.88,-1.08)**@{-},
 (0.3,-.35)*+{y},
(4.0,3.3)*+{(\frac1{m-1},(-1)^{m-1})},
(-3.7,3.3)*+{(\frac1{m+1},(-1)^{m+1})},
(0,-3.5)*+{(\frac1m,(-1)^m)},
(-1.2,-1)*\cir<8mm>{}="*",
(-1.2,-1);
(-1.4,-1.12)**@{-},
(0,-3);
(-.2,-3.12)**@{-},
(-1.4,-1.12);
(-.2,-3.12)**\dir{.},?<(0)*\dir{<},?>(1)*\dir{>},
(-1.3,-2.3)*+{>t},
 \endxy
\begin{fig} \label{f7}
\end{fig}

Let $y=(y_1, y_2)$ be the point with $y_1 = \frac{1}{m}$,  so that the segment $[x,y]$
with endpoints $x$ and $y$ is perpendicular to $I_m$.
If $t \geq \frac{r}{3}$ then, by \eqref{m}, (since $m>3$),
$$
d=d(x,y)\geq t \cdot \sin \be_m \geq \frac{r}{3} \cdot \frac{1}{3m(m+1)} \geq
\frac{r}{3} \cdot \frac {1}{3} \cdot \frac {1}{4}\cdot \frac{r}{2} \cdot \frac{1}{5} \cdot
\frac{r}{2}= \frac{r^3}{720}.
$$

If $d<r$ then $y \in B (x, r)$ and $f(y)=d \geq \frac{r^3}{720}$.  If $d \geq r$,  let
$y_r$ denote a point in the segment with endpoints $x$ and $y$ so that $d(x, y_r)=r$.
Then $f(y_r)=r$ and $f(B(x, r))\supset f ([x, y_r))\supset [0, r)$.

Next we consider the case when $t< \frac{r}{3}$, as illustrated on Figure~\ref{f8}

\hspace{27mm}  \xy <0.7cm,0cm>:
 (-3.6,3) ;
 (0,-3)**@{-},
 (0,-3) ;
 (3.8,3)**@{-},
 (0,3) ;
 (0,-3)**@{.},
 (-1.2,-1)*+\dir{*},
 (-1.2,-1);
 (0,-1)**@{-},
 (-1.5,-1)*+{x},
 (-.7,-.9)*+{\eta},
 (0,1.38)*+\dir{*},
 (0,1.38);
 (-1.2,-1)**@{-},
 (0,1.38);
 (-1.93,.22)**@{-},
 (-1.2,1.2)*+{f(z)},
 (0.3,1.4)*+{z},
 (0.3,0)*+{\gamma},
 (-.9,0.2)*+{r},
(4.0,3.3)*+{(\frac1{m-1},(-1)^{m-1})},
(-3.7,3.3)*+{(\frac1{m+1},(-1)^{m+1})},
(0,-3.5)*+{(\frac1m,(-1)^m)},
(-1.2,-1);
(-1.4,-1.12)**@{-},
(0,-3);
(-.2,-3.12)**@{-},
(-1.4,-1.12);
(-.2,-3.12)**\dir{.},?<(0)*\dir{<},?>(1)*\dir{>},
(-1.3,-2.3)*+{>t},
 \endxy
\begin{fig} \label{f8}
\end{fig}

Then
$$
\eta = d(x, \{v=(v_1,v_2):v_1 = \frac{1}{m}\})
\leq 2t \sin \be_m \leq 2 \cdot \frac{r}{3} \cdot \frac{1}{m(m+1)}
\leq \frac{r}{3}.
$$
Hence
$$
\gamma = \sqrt{r^2- \eta^2} \geq \sqrt{r^2 - \frac{r^2}{9}} =r \frac{\sqrt{8}}{3}.
$$
Thus there exists $z \in B(x, r)$ with $z_1 = \frac{1}{m}$ and
$z_2 \geq x_2 + r \frac{\sqrt{8}}{3}\geq -1 + \frac{2}{3}r$.
We have
$$
f(z)= d(z,K) \geq \frac{2}{3}r \cdot \sin \be_m \geq \frac{2}{3}r \cdot
\frac{1}{3m(m+1)} \geq \frac{2}{3}r \cdot \frac{1}{3} \cdot \frac{1}{4}\cdot  \frac{r}{2}
\cdot \frac{1}{5} \cdot\frac{r}{2} = \frac{r^3}{360}.
$$
Thus for all $t \in [0,1]$ we conclude that $f(B(x,r))\supset [0, \frac{r^3}{720}]$.

A similar computation shows that $f(B(x,r))$ contains also a sufficiently  large
negative interval, so that \eqref{rcube} holds.
\end{proof}

\begin{step} \lb{2}
If $d(x, K)= d>0$ then
\begin{itemize}
\item[(a)]
if $f(x)>0$ then
\begin{equation}\lb{xclose}
f(B(x, r))\supset (\max(d-r,0), d +\frac{r^3}{480});
\end{equation}
\item[(b)]
if $f(x)< 0$ then
\begin{equation}\lb{xclosen}
f(B(x, r))\supset (-d-\frac{r^3}{480}, \min(d+r,0)).
\end{equation}
\end{itemize}
\end{step}

\begin{proof}[Proof of Step~\ref{2}]
We will assume without loss of generality that $f(x)>0$.  The case when  $f(x)<0$ is
proven identically. To prove \eqref{xclose} we will consider two cases.  First we
assume that $x= (x_1, x_2)$ where $|x_1|,\ |x_2|\leq 1$.  Since $x_1 \neq 0$, this
implies that there exists $m \in \bbZ$ (say, $m>0$) so that $x$ lies inside the
triangle with vertices $(\frac{1}{m-1},(-1)^{m-1})$,
$(\frac{1}{m},(-1)^m)$, $(\frac{1}{m+1},(-1)^{m+1})$, as illustrated on Figure~\ref{f3}


\hspace{4mm}
\xy <0.6cm,0cm>:
 (-3.6,3) ;
 (0,-3)**@{-},
 (10,3) ;
 (14,-3)**@{-},
 (14,-3) ;
 (18,3)**@{-},
 (14,3) ;
 (14,-3)**@{.},
 (0,-3) ;
 (3.8,3)**@{-},
 (0,-1) ;
 (0,-3)**@{.},
 (-1.2,-1)*+\dir{*},
 (-1.6,-1)*+{z},
(4.0,3.4)*+{(\frac1{m-1},(-1)^{m-1})},
(-3.7,3.4)*+{(\frac1{m+1},(-1)^{m+1})},
(0,-3.5)*+{(\frac1m,(-1)^m)},
(14,-3.5)*+{(\frac1n,(-1)^n)},
(-1.02,-1.3);
(-.72,-1.12)**@{-},
(-.9,-.82);
 (-.72,-1.12)**@{-},
 (-1.2,-1);
(-.2,-.4)**@{-},
(-1.26,1.36);
(-.2,-.4)**@{-},
(-1.26,1.8)*+{d'},
(-0.2,-.4)*+\dir{*},
(0,-.7)*+{x},
(-.6,-2);(0.63,-2)**\crv{(0,-1.65)},
(-.8,.6);(-.2,.77)**\crv{(-.5,.8)},
(-2,-2.5);(-.2,-2.2)**\dir{-} ?>(1)*\dir{>},
(-2.4,-2.5)*+{\beta_m},
(2,-2.5);(.2,-2.2)**\dir{-} ?>(1)*\dir{>},
(2.4,-2.5)*+{\alpha_m},
(-2.5,-.5);(-.5,.4)**\dir{-} ?>(1)*\dir{>},
(-2.7,-.6)*+{\beta_m},
(-0.2,2)*+\dir{*},
(14,3.4)*+{v},
(0.1,2.4)*+{y},
(0.2,.8)*+{r},
(-.2,2);
(-.2,-.4)**@{-},
(-.2,2);
(-2.26,.76)**@{-},
(14,-.4);
(-.2,-.4)**@{.},
 \endxy
\begin{fig} \label{f3}
\end{fig}

If $\frac{1}{m+1}\geq \frac{r}{4}$ we consider a point $y=(y_1,y_2) \in
\overline{B(x,r)}$ so that $y_1 = x_1,\ y_2=x_2+r$.   Then
\begin{equation*}
\begin{split}
f(y) &= d' = d(y, K) \geq d + r\sin \be_m\\
&\geq d + r \frac{1}{3m(m+1)} \geq d + r \cdot\frac{1}{3}\cdot \frac{r}{4} \cdot \frac{r}{4}\\
&= d + \frac{r^3}{48}.
\end{split}
\end{equation*}

If $\frac{1}{m+1} < \frac{r}{4}$ then we proceed similarly to Steps~\ref{3} and
\ref{4}.  Let $n$ be the smallest odd number so that
$$
\frac{1}{n-1}> \frac{r}{2}
$$
Since $x_1 > 0$, we see that there exists $v=(v_1, v_2) \in B(x,r)$ so that
$v_1=\frac{1}{n}, v_2 \geq x_2 + \frac{r}{2}$.  Then, as before,
$$
f(v)=d(v,K)\geq \frac{r}{2}\cdot \sin \be_n \cdot\frac{r}{2} \cdot\frac{1}{3n(n+1)} \geq
\frac{r^3}{480}.
$$

Now let $z$ be that point on the interval $I_n$ so that $d(x, z)=d$.  If $r>d$  then
$z \in B(x,r)$ and thus $[0,d] \subset f (B(x, r))$.  If $r \leq d$ then $B(x, r)$
contains a subinterval of length r of the interval $[x, z]$ and $f(B(x, r)) \supset
(d-r,d]$.
Thus \eqref{xclose} is satisfied.

Next we consider the case when $|x_2| > 1$, as illustrated on Figure~\ref{f4}.


 \hspace{30mm} \xy <1cm,0cm>:
(0,0) ;
 (1,2)**@{-},
 (1,2) ;
 (2,0)**@{-},
 (2,0) ;
 (3.2,2)**@{-},
 (3.2,2);
 (4.5,0)**@{-},
 (6.5,0) ;
 (4.5,0)**@{-},
(2.5,4)*+\dir{*},
(2.5,4); (1,2)**@{-},
(2.5,4)*\cir<12mm>{},
(2.5,4); (2.5,5.1)**@{-},
(1.4,4); (2.5,5.1)**@{--},
(1,2); (2.5,5.1)**@{-},
(2.5,4); (1.4,4)**@{-},
(2.7,4.6)*+{r},
(2.7,3.9)*+{x},
(1.2,4)*+{z},
(1.4,3.3)*+{v},
(2.6,5.3)*+{y},
(2.5,4);(1.6,3.3)**@{--},
(2,3)*+{d},
\endxy

\vspace{1mm}
\begin{fig}  \label{f4}
\end{fig}

Then, as above, it is clear that $f(B(x, r)) \supset (\max(0, d-r),d]$.  Further
there exists $y=(y_1,y_2)\in \overline{B(x, r)}$ so that $y_1 = x_1,\ y_2 = x_2{+r}$.
Let $z = (z_1, z_2)\in \overline{B(x, r)}$ be such that $z_1 = x_1{-r},\ z_2 = x_2$,
and $v \in \overline{B(x, r)}$ be so that $d(v, x)=r$ and $v$ lies on the shortest
path from $y$ to $K$.  Then
$$
f(y)= d(y, K)=d(v, K)+d(v, y) \geq(d-r)+ d(z, y) = d-r+ \sqrt{2}r \geq d + \frac{r}{3}.
$$
Thus \eqref{xclose} holds.

The case when $|x_1|>1$ follows from very similar considerations, which ends the
proof of Step~\ref{2}.
\end{proof}

As an immediate corollary of Step~\ref{2} we obtain the following:
\begin{step} \lb{2a}
If $d(x, K)= d>0$ and $r \leq \min (\frac{1}{10}, d)$ then
$$
f(B(x, r))\supset B (f(x),\frac{r^3}{480}).
$$
\end{step}

\begin{step}\lb{5}
If $d(x, K)= d>0$ and $\frac{1}{10}\geq r>d$ then
\begin{equation} \lb{last}
f(B(x, r)) \supset B(f(x),\frac{r^3}{16000}).
\end{equation}
\end{step}

\begin{proof}[Proof of Step~\ref{5}]
We start from the trivial observation that when
$\frac{1}{10} \geq r > d$ then, by Step~\ref{2a},
$$
f(B(x, r))\supset f (B(x, d))\supset B(f(x),\frac{d^3}{480}).
$$
Thus, if
\begin{equation} \lb{rsmall}
\frac{r^3}{16000} \leq \frac{d^3}{480},
\end{equation}
then \eqref{last} is satisfied.  Equation \eqref{rsmall} is true when $r \leq 2d$.
Thus, next we assume that
$$
\frac{1}{10} \geq r \geq 2d.
$$
We will also assume, without loss of generality, that $f(x)>0$.
Let $y \in K$ be such that $d(x, y)=d$.  Then $B(x, r)\supset B(y, r-d)$.  By
Steps~\ref{3} and \ref{4}, \eqref{rcube} holds and we have
$$
f(B(y, r-d)) \supset (-\frac{(r-d)^3}{2000},0].
$$
Note that
$$
\frac{(r-d)^3}{2000} \geq \frac{(\frac{1}{2}r)^3}{2000}=\frac{r^3}{16000}.
$$
Thus
$$
f(B(x, r))\supset (-\frac{r^3}{16000}, d].
$$
On the other hand, by Step~\ref{2},
$$
f(B(x, r)) \supset (0, d + \frac{r^3}{480}).
$$
Thus \eqref{last} is satisfied.\end{proof}

This ends the proof that $f$ is co-uniformly continuous with the modulus $\omega$
defined in \eqref{omega}.
\end{proof}


\def\polhk#1{\setbox0=\hbox{#1}{\ooalign{\hidewidth
  \lower1.5ex\hbox{`}\hidewidth\crcr\unhbox0}}} \def\cprime{$'$}

\end{document}